\documentstyle[11pt,bezier]{article}
\catcode`\@=11
\ifcase\@ptsize
 \font\tenmsa=msam10
 \font\sevenmsa=msam7
 \font\fivemsa=msam5
 \font\tenmsb=msbm10
 \font\sevenmsb=msbm7
 \font\fivemsb=msbm5
 \font\teneu=eufm10
 \font\seveneu=eufm7
 \font\fiveeu=eufm5
 \font\tenib=cmmib10
 \font\sevenib=cmmib7
 \font\fiveib=cmmib5
\or
 \font\tenmsa=msam10 scaled \magstephalf
 \font\sevenmsa=msam7 scaled \magstephalf
 \font\fivemsa=msam5 scaled \magstephalf
 \font\tenmsb=msbm10 scaled \magstephalf
 \font\sevenmsb=msbm7 scaled \magstephalf
 \font\fivemsb=msbm5  scaled \magstephalf
 \font\teneu=eufm10  scaled \magstephalf
 \font\seveneu=eufm7  scaled \magstephalf
 \font\fiveeu=eufm5   scaled \magstephalf
 \font\tenib=cmmib10  scaled \magstephalf
 \font\sevenib=cmmib7  scaled \magstephalf
 \font\fiveib=cmmib5   scaled \magstephalf
\or
 \font\tenmsa=msam10 scaled \magstep1
 \font\sevenmsa=msam7 scaled \magstep1
 \font\fivemsa=msam5  scaled \magstep1
 \font\tenmsb=msbm10 scaled \magstep1
 \font\sevenmsb=msbm7 scaled \magstep1
 \font\fivemsb=msbm5  scaled \magstep1
 \font\teneu=eufm10   scaled \magstep1
 \font\seveneu=eufm7 scaled \magstep1
 \font\fiveeu=eufm5 scaled \magstep1
 \font\tenib=cmmib10     scaled \magstep1
 \font\sevenib=cmmib7   scaled \magstep1
 \font\fiveib=cmmib5   scaled \magstep1
\fi
\newfam\msafam
\textfont\msafam=\tenmsa  \scriptfont\msafam=\sevenmsa
  \scriptscriptfont\msafam=\fivemsa
\newfam\msbfam
\textfont\msbfam=\tenmsb  \scriptfont\msbfam=\sevenmsb
  \scriptscriptfont\msbfam=\fivemsb
\def\Bbb{\ifmmode\let\next\Bbb@\else
 \def\next{\errmessage{Use \string\Bbb\space only in math mode}}\fi\next}
\def\Bbb@#1{{\Bbb@@{#1}}}
\def\Bbb@@#1{\fam\msbfam#1}
\newfam\eufam
\textfont\eufam=\teneu  \scriptfont\eufam=\seveneu
  \scriptscriptfont\eufam=\fiveeu
\def\frak{\ifmmode\let\next\frak@\else
 \def\next{\errmessage{Use \string\frak\space only in math mode}}\fi\next}
\def\frak@#1{{\frak@@{#1}}}
\def\frak@@#1{\fam\eufam#1}
\newfam\ibfam
\textfont\ibfam=\tenib  \scriptfont\ibfam=\sevenib
  \scriptscriptfont\ibfam=\fiveib
\def\bold{\ifmmode\let\next\bold@\else
 \def\next{\errmessage{Use \string\bold\space only in math mode}}\fi\next}
\def\bold@#1{{\bold@@{#1}}}
\def\bold@@#1{\fam\ibfam#1}
\def\hexnumber@#1{\ifcase#1 0\or 1\or 2\or 3\or 4\or 5\or 6\or 7\or 8\or
 9\or A\or B\or C\or D\or E\or F\fi}
\def\newsymbolb#1#2#3#4{\mathchardef#1="#2\hexnumber@\msbfam#3#4}
\def\newsymbola#1#2#3#4{\mathchardef#1="#2\hexnumber@\msafam#3#4}
\newsymbola{\leqslant}336
\newsymbola{\geqslant}33E
\newsymbola{\square}003
\newsymbola{\lessgtr}337
\newsymbola{\gtrless}33F
 \let\leq\leqslant
 \let\geq\geqslant
\font\fraksect=eufm10 scaled 1728
\font\fraknote=eufm8
\font\frakssect=eufm10 scaled 1440
\def\hybrid{\topmargin 0pt      \oddsidemargin 0pt
        \headheight 0pt \headsep 0pt
        \textwidth 160true mm       % US paper
        \textheight 231true mm         % US paper
        \marginparwidth 0.0in
        \parskip 0pt plus 1pt   \jot = 1.5ex}
\catcode`\@=11
\def\marginnote#1{}

\newcount\hour
\newcount\minute
\newtoks\amorpm
\hour=\time\divide\hour by60
\minute=\time{\multiply\hour by60 \global\advance\minute by-\hour}
\edef\standardtime{{\ifnum\hour<12 \global\amorpm={am}%
        \else\global\amorpm={pm}\advance\hour by-12 \fi
        \ifnum\hour=0 \hour=12 \fi
        \number\hour:\ifnum\minute<10 0\fi\number\minute\the\amorpm}}
\edef\militarytime{\number\hour:\ifnum\minute<10 0\fi\number\minute}

\def\draftlabel#1{{\@bsphack\if@filesw {\let\thepage\relax
   \xdef\@gtempa{\write\@auxout{\string
      \newlabel{#1}{{\@currentlabel}{\thepage}}}}}\@gtempa
   \if@nobreak \ifvmode\nobreak\fi\fi\fi\@esphack}
        \gdef\@eqnlabel{#1}}
\def\@eqnlabel{}
\def\@vacuum{}
\def\draftmarginnote#1{\marginpar{\raggedright\scriptsize\tt#1}}

\def\draft{\oddsidemargin -.5truein
        \def\@oddfoot{\sl preliminary draft \hfil
        \rm\thepage\hfil\sl\today\quad\militarytime}
        \let\@evenfoot\@oddfoot \overfullrule 3pt
        \let\label=\draftlabel
        \let\marginnote=\draftmarginnote
   \def\@eqnnum{(\theequation)\rlap{\kern\marginparsep\tt\@eqnlabel}%
\global\let\@eqnlabel\@vacuum}  }
\newcounter{app}
\newcounter{sapp}[app]

\makeatletter
\newdimen\normalarrayskip              % skip between lines
\newdimen\minarrayskip                 % minimal skip between lines
\normalarrayskip\baselineskip
\minarrayskip\jot
\newif\ifold             \oldtrue            
\def\arraymode{\ifold\relax\else\displaystyle\fi} % mode of array entries
\def\eqnumphantom{\phantom{(\theequation)}}     % right phantom in eqnarray
\def\@arrayskip{\ifold\baselineskip\z@\lineskip\z@
     \else
     \baselineskip\minarrayskip\lineskip2\minarrayskip\fi}
\def\@arrayclassz{\ifcase \@lastchclass \@acolampacol \or
\@ampacol \or \or \or \@addamp \or
   \@acolampacol \or \@firstampfalse \@acol \fi
\edef\@preamble{\@preamble
  \ifcase \@chnum
     \hfil$\relax\arraymode\@sharp$\hfil
     \or $\relax\arraymode\@sharp$\hfil
     \or \hfil$\relax\arraymode\@sharp$\fi}}
\def\@array[#1]#2{\setbox\@arstrutbox=\hbox{\vrule
     height\arraystretch \ht\strutbox
     depth\arraystretch \dp\strutbox
     width\z@}\@mkpream{#2}\edef\@preamble{\halign \noexpand\@halignto
\bgroup \tabskip\z@ \@arstrut \@preamble \tabskip\z@ \cr}%
\let\@startpbox\@@startpbox \let\@endpbox\@@endpbox
  \if #1t\vtop \else \if#1b\vbox \else \vcenter \fi\fi
  \bgroup \let\par\relax
  \let\@sharp##\let\protect\relax
  \@arrayskip\@preamble}
\def\eqnarray{\stepcounter{equation}%
              \let\@currentlabel=\theequation
              \global\@eqnswtrue
              \global\@eqcnt\z@
              \tabskip\@centering
              \let\\=\@eqncr
              $$%
 \halign to \displaywidth\bgroup
    \eqnumphantom\@eqnsel\hskip\@centering
    $\displaystyle \tabskip\z@ {##}$%
    &\global\@eqcnt\@ne \hskip 2\arraycolsep
         $\displaystyle\arraymode{##}$\hfil
    &\global\@eqcnt\tw@ \hskip 2\arraycolsep
         $\displaystyle\tabskip\z@{##}$\hfil
         \tabskip\@centering
    &{##}\tabskip\z@\cr}
\makeatother
\relax

\hybrid

\begin{document}
\def\bea{\begin{eqnarray}}
\def\eea{\end{eqnarray}}
\def\beq{\begin{equation}}          \def\bn{\beq}
\def\eeq{\end{equation}}            \def\ed{\eeq}
\def\nn{\nonumber}                  \def\g{\gamma}
\def\Uq{U_q(\widehat{\frak{sl}}_2)}
\def\Uqp{U_q(\widehat{\frak{sl}}'_2)}
\def\Uqd{U^{*}_q(\widehat{\frak{sl}}_2)}
\def\uq{U_q({sl}_2)}
\def\uqd{U^*_q({sl}_2)}
\def\slaff{\frak{sl}^\prime_2}
\def\aff{\widehat{\frak{sl}}_2}
\def\ot{\otimes}
\def\sk#1{\left({#1}\right)}
\def\id{\mbox{\rm id}}
\def\tr{\mbox{\rm tr}}
\def\tah{\mbox{\rm th}}
\def\sh{\mbox{\rm sh}}
\def\ch{\mbox{\rm ch}}
\def\ctg{\mbox{\rm ctg}}
\def\cth{\mbox{\rm cth}}
\def\tg{\mbox{\rm tg}}
\def\th{\mbox{\rm th}}
\def\sign{\mbox{\rm sign}}
\def\qdet{\mbox{\rm q-det}}
\def\Re{{\rm Re}\,}
\def\Im{{\rm Im}\,}
\def\RR{\Bbb{R}}
\def\ZZ{\Bbb{Z}}
\def\CC{\Bbb{C}}
\def\r#1{\mbox{(}\ref{#1}\mbox{)}}
\def\d{\delta}
\def\D{\Delta}
\def\da{{\partial_\alpha}}
\let\da=p
\def\Ps{\Psi^{*}}
\def\R{{\cal R}}
\def\Ga#1{\Gamma\left(#1\right)}
\def\si#1{\sin\pi\left(#1\right)}
\def\ex#1{\exp\left(#1\right)}
\def\ep{\varepsilon}
\def\eps{\epsilon}
\def\ve{\ep}
\def\fract#1#2{{\mbox{\footnotesize $#1$}\over\mbox{\footnotesize $#2$}}}
\def\stackreb#1#2{\ \mathrel{\mathop{#1}\limits_{#2}}}
\def\res#1{\stackreb{\mbox{\rm res}}{#1}}
\def\lim#1{\stackreb{\mbox{\rm lim}}{#1}}
\def\Res#1{\stackreb{\mbox{\rm Res}}{#1}}
\let\dis=\displaystyle
\def\ee{{\rm e}}
\def\D{\Delta}
\def\HH{{\cal H}}
\renewcommand{\theequation}{{\thesection}.{\arabic{equation}}}
\def\Y-{\widehat{Y}^-}
\font\fraksect=eufm10 scaled 1728
\font\fraknote=eufm8
\font\frakssect=eufm10 scaled 1440
\def\DYsect{\widehat{DY(\hbox{\fraksect sl}_2)}}
\def\DYssect{\widehat{DY(\hbox{\frakssect sl}_2)}}
\def\Aya{{\cal A}_{\tih}(\widehat{\frak{sl}_2})}
\def\Ayap{{\cal A}^+_{\tih}(\widehat{\frak{sl}_2})}
\def\Ayam{{\cal A}^-_{\tih}(\widehat{\frak{sl}_2})}
\def\Ayasect{{\cal A}_{\tih}(\widehat{\hbox{\fraksect sl}_2})}
\def\Ayassect{{\cal A}_{\tih}(\widehat{\hbox{\frakssect sl}_2})}
\def\Ael{{\cal A}_{\tih,\eta}(\widehat{\frak{sl}_2})}
\def\Aelx{{\cal A}_{\tih,\xi}(\widehat{\frak{sl}_2})}
\def\Aelsect{{\cal A}_{\tih,\eta}(\widehat{\hbox{\fraksect sl}_2})}
\def\Aelnote{{\cal A}_{\tih,\eta}(\widehat{\hbox{\fraknote sl}_2})}
\def\Aelssect{{\cal A}_{\tih,\eta}(\widehat{\hbox{\frakssect sl}_2})}
\def\Apq{{\cal A}_{q,p}(\widehat{\frak{sl}_2})}
\def\Apqsect{{\cal A}_{q,p}(\widehat{\hbox{\fraksect sl}_2})}
\def\DY{\widehat{DY(\frak{sl}_2)}}
\def\Yd{\DY}
\def\Ydd{\DY}
\let\z=z
\let\b=z
\def\u{{u}}
\def\v{{v}}
\def\g{\gamma}
\def\la{\lambda}
\let\hsp=\qquad
\def\he{{\hat e}}
\def\hf{{\hat f}}
\def\hh{{\hat t}}
\def\ha{{\hat a}}
\def\hb{{\hat b}}
\def\hk{{\hat t}}
\def\ka{\kappa}
\def\kk{{\hat \ka}}
\def\hkp{{\hat t}}
\def\hhh{{\hat h}}
\def\vvv{\overline{\varphi}}
\def\vac{|\mbox{vac}\rangle}
\def\lvac{\langle\mbox{vac}|}
\def\vpint{-\!\!\!\!\!\!\!\int_{-\infty}^\infty}
\def\vpinto{-\!\!\!\!\!\!\!\int_{0}^\infty}
\def\vpintss{{\vpint\cdots\ \vpint}}
\def\vpintoo{{\vpinto\cdots\ \vpinto}}
\def\intt{\int_{-\infty}^\infty}
\def\vpints{-\!\!\!\!\!\!\!\int_{-\infty}^{(\la-\mu)/2}}
\def\la{\lambda}
\def\tih{\hbar}
\def\h{\hbar}
\def\Ev{{\cal E}v}
\def\FDY{F\left[\DY\right]}
\def\FDYsect{F\left[\DYsect\right]}
\def\stackupb#1#2#3{\ \mathrel{\mathop{#1}\limits_{#2}^{#3}}}
\def\feq#1#2{\stackupb{\ravnodots}{#1}{#2}}
\def\vac{|\mbox{vac}\rangle}
\def\lvac{\langle\mbox{vac}|}
\def\cint{\int_\infty^{0+}}
\def\nint{\int^{+\infty}_{0}}
\def\mint{\int_{-\infty}^{0}}
\let\rvac=\vac
\def\Cup{\bigcup\limits}
\def\yh{\nu}
\def\E{E}
\def\F{F}
\def\H{H}
\hyphenation{co-mul-tip-li-ca-tion auto-mor-phism
rep-re-sen-ta-tion
rep-re-sen-ta-tions de-ge-ne-ra-tion sub-al-geb-ras
com-mu-ta-tion con-ve-ni-ent anti-iso-mor-phism ge-ne-ra-tors
cor-res-pon-ding for-ma-lism co-or-di-na-tes com-pac-ti-fied
ana-ly-ti-cal fac-to-ri-za-tion}
%%%%%%%%%%%%%%%%%%%%%%%%%%%%%%%%%%%%%%%%%%%%%%%%%%%%%%%%%%%%%%%%%%%%%%
%%%%%%%%%%%%%%%%%   End of Personal difinitions    %%%%%%%%%%%%%%%%%%%
%%%%%%%%%%%%%%%%%%%%%%%%%%%%%%%%%%%%%%%%%%%%%%%%%%%%%%%%%%%%%%%%%%%%%%

\makeatletter

\def\Cup{\bigcup\limits} \def\Sum{\sum\limits} \def\Prod{\prod\limits}
\let\Int\int \def\int{\Int\limits}
\let\Oint\oint \def\oint{\Oint\limits}
\def\sltwo{\frak{sl}_2}
\def\refskips{\parskip 0pt plus .1pt\partopsep0pt}
\def\vsk#1>{\vskip#1\baselineskip}
\let\vp\vphantom \let\hp\hphantom
\def\ftext#1{{\let\thefootnote\relax\footnotetext{\vsk-.8>\noindent #1}}}
\parindent 16pt
\def\@oddfoot{\footnotesize\hfil\rm\thepage\hfil}
\let\@evenfoot\@oddfoot
\let\tsize\textstyle
\def\qed{\hbox{}\nobreak\hfill\nobreak{\m@th$\,\square$}}
\def\Ucup{{\tsize\bigcup}}
\def\Zln{{\cal Z}^\ell_n}
\def\Zlnb{\rlap{$\,\overline{\!\phantom{\cal Z}}$}\Zln}
\def\lg{\frak l}
\def\mg{\frak m}

\def\@xthm#1#2{\@begintheorem{#2}{\csname the#1\endcsname.}\ignorespaces}
\def\@ythm#1#2[#3]{\@opargbegintheorem{#2}{\csname the#1\endcsname}{#3.}%
 \ignorespaces}

\makeatother

\newenvironment{abst}{\begingroup\narrower\noindent{\bf Abstract.}\enspace
 \ignorespaces}{\endgraf\endgroup}
\newenvironment{proof}{\begin{trivlist}\item{\it Proof.}\enspace\ignorespaces}
 {\qed\end{trivlist}}
\newenvironment{remark}{\begin{trivlist}\item{\sl Remark.}\enspace
 \ignorespaces}{\end{trivlist}}
\newenvironment{comment}{\begin{trivlist}\item\raggedright\tt}{\end{trivlist}}

\newtheorem{prop}{Proposition}
\newtheorem{lemma}[prop]{Lemma}
\newtheorem{corol}[prop]{Corollary}

\newcounter{parag}
\def\razdel{\addtocounter{parag}{1}\thesection.\theparag. }

%\begin{titlepage}
\begin{center}
%{Revised \today}\\
q-alg/9712057
\hfill ITEP-TH-66/97\\
\hfill Roma1-1189/97\\
%\hfill q-alg/9712057\\
\bigskip\bigskip
{\Large\bf Yangian Algebras and Classical Riemann Problems}\\
\bigskip
\bigskip
{\large S. Khoroshkin$\,^\star$\footnote{E-mail: khoroshkin@vitep1.itep.ru},
D. Lebedev$\,^\star$\footnote{E-mail: lebedev@vitep1.itep.ru} and
S. Pakuliak$\,^{\star\,*\,\diamond}$\footnote{E-mail:
pakuliak@thsun1.jinr.ru}}
\bigskip
\end{center}
\bigskip

$\quad^\star${\it Institute of Theoretical \& Experimental Physics,
117259 Moscow, Russia}

$\quad^*${\it Bogoliubov Laboratory of Theoretical Physics, JINR,
141980 Dubna, Russia}

$\quad^\diamond${\it Bogoliubov Institute of Theoretical Physics,
252143, Kiev, Ukraine}
\bigskip
\bigskip

\begin{abstract}
%The algebra of $L$-operators corresponding to the rational sulution
%f quantum Yang-Baxter equation is considered.
%tarting from the Yang solution of quantum Yang--Baxter equation
%%
 We investigate different
Hopf algebras associated to Yang's solution of quantum
 Yang--Baxter equation.
It is shown that for the
precise definition of the algebra
one needs the commutation relations for the deformed algebra
 of formal currents and the specialization of the
 Riemann problem for the currents.
Two different Riemann
problems are considered.
They lead to the central extended
Yangian double associated with
$\frak{sl}_2$  and to the degeneration of scaling limit of elliptic
affine algebra.
Unless the defining relations for the generating functions of
the both algebras coincide their properties and  the theory
of infinite-dimensional representations are quite different.
%ue to different analytical properties and meaning of the
%enerating functions.
 We discuss also
the Riemann problem for
twisted algebras and for
scaled elliptic algebra.
\end{abstract}
%\end{titlepage}
%\clearpage

\setcounter{equation}{0}
\setcounter{footnote}{0}

\section{Introduction}
The Yangian $Y(\frak{g})$, where $\frak{g}$ is a simple Lie algebra
 was introduced by V.~Drinfeld \cite{Dr2} as a Hopf algebra  such that
a quantization of the Yang rational solution of classical Yang-Baxter
 equation can be done in tensor category of finite-dimensional
 representations of $Y(\frak{g})$.
  Later the algebraical structure of the quantum double
$DY(\frak{g})$ of the Yangian was
 studied in \cite{KT} and  the corresponding universal $R$-matrix was
 calculated explicitly. The Yangian double admits a central extension
$\widehat{DY(\frak{g})}$ \cite{K} (see also \cite{IK}) and  the intertwining
 operators for its infinite-dimensional representations
can be used
 for the calculation of form-factors of  local operators in $SU(2)$
 invariant Thirring model \cite{KLP4}.

Recently, the detail analysis of the scaling limit $\Ael$ of the elliptic
algebra was done in \cite{KLP3}. The corresponding intertwining operators
 can be used for the calculations
of  the form-factors in $XXZ$ model in the gapless regime \cite{JM}
and in the
 Sine-Gordon model \cite{L}. It happens that rational degeneration
 $\Aya$ of the
 algebra $\Ael$ ($\eta\to0$)
can be described on the level of generating functions
 ($L$-operators) by the same set of the relations as central extended
 Yangian double $\DY$ while the structure of the algebras and their
 infinite-dimensional representations look very different.

For an investigation of this phenomena we go back to the original
 ideology of classical inverse scattering method, where the Riemann
 problem of factorizing the matrix valued function into the product
 of functions analytical in certain domains plays the crucial role
\cite{FTbook}.
We claim that the complete structure of the quantum algebra, including
 its  coalgebraic structure, can be given by the following data:
 a formal algebra of (deformed) currents and a specific Riemann
 problem. For the algebras related to the
$\frak{sl}_2$ case which we study here it means that
 they are given by the formal relations for
the total currents $e(u)$, $f(u)$ and $h^\pm(u)$
%%%%%%%%%%%%%%%    (e,f,h are generators of $\frak{sl_2}$,
(here $u$ is a spectral parameter)
 and the decompositions $e(u)=e^+(u)-e^-(u)$, $f(u)=f^+(u)-f^-(u)$.
The operator valued generating functions
of the spectral parameter
$e^\pm(u)$, $f^\pm(u)$ are given as certain integral transforms
 of the total currents $e(u)$ and $f(u)$
which can be  uniquely  determined by the
 predicted analytical properties of $e^\pm(u)$ and $f^\pm(u)$ (see Section 3).
 The specification of the Riemann problem uniquely defines the
 precise expressions of the currents as generating functions of the
 elements of the algebra. It converts the relation between the
 currents into the relations between the generators of the algebra.
 Moreover it defines the comultiplication structure for the
 quantum algebra. We assume here that $e^\pm(u)$, $f^\pm(u)$ and
 $h^\pm(u)$     are Gauss coordinates of some $L$-operators
and  use the  universal comultiplication formulas for
Gauss coordinates (see Section 7).

For both algebras $\DY$ and $\Aya$ the commutation relations between
 the currents can be computed by a standard procedure called
Ding-Frenkel isomorphism \cite{DF}. Actually the calculations are
 formal algebraical manipulations which use only the structure
 of the Yang (quantum) $R$-matrix $R(u)=1+\h P/u$ where $P$ is a flip
 and thus coincide. This is done in Section 2.

The Riemann problem for $\DY$ means a decomposition of a function with finite
 number of singularities into the sum of functions analytical at zero and
 at infinity. The solution is given by Cauchy integrals
 $e^\pm(u)=\oint\frac{e(v)dv}{2\pi i(u-v)}$, where the closed contour
including zero goes from the left or from the right of the point $u$. The
elements of the corresponding algebra are Taylor coefficients of $e^\pm(u)$,
$f^\pm(u)$ and of $h^\pm(u)$. They generate precisely $\DY$.  The Riemann
problem which corresponds to the algebra $\Aya$ means a decomposition of a
 function vanishing at infinity into the sum of functions analytical in
 some half-planes. The solution is given by the Cauchy integral with the
 same kernel over the contour being the straight line parallel to the real
 axis going lower or upper $u$. The generators of the algebra, indexed by
 real numbers are the coefficients of inverse Laplace transforms of the
  currents  $e^\pm(u)$, $f^\pm(u)$ and of $h^\pm(u)$ (Section 3).
 The analysis of the natural completions of the two algebras and of their
  properties show that they cannot be transformed one into another
  by means of projective transforms as well as corresponding Riemann
  problems are essentially different due to different asymptotical
  conditions for the involved functions (Section 4).

These two Riemann
  problems may be modified by imposing other asymptotical behaviors
  for $e^\pm(u)$, $f^\pm(u)$. Thus we get the twisted variants of these
  two algebras and as a limit the popular algebras with a simple
  comultiplication introduced by Drinfeld (new realization) (see Section 7).
 We demonstrate also in this Section that further trigonometric generalization
 of the Yangian current algebra essentially leads to the algebra defined in
\cite{KLP3} and an application of the Riemann problem for a strip
 produces the  scaled elliptic algebra $\Ael$.

  It is worth to mention that the algebra $\Aya$ fits much more for the
  applications in quantum integrable field theories. The main advantage
  is that the algebra $\Aya$ is graded while $\DY$ is a filtered algebra.
  Moreover the grading operator $d$ can be diagonalized in
  infinite-dimensional representations of $\Aya$ such that its trace has a
  sense contrary to the case of $\DY$ where only the ratio of the traces
  is well defined (see Section 6). The application of the universal
  $R$-matrix to finite-dimensional representations gives integral forms of
   the corresponding $R$-matrices (Section 5). One may consider the universal
 $R$-matrix for $\Aya$ as another quantization of classical rational solution
  of the Yang-Baxter equation.

\setcounter{equation}{0}
\setcounter{parag}{0}
\section{Algebra of $L$-Operators}

{\bf \razdel}
The aim of this section is to develop well known technique of the
$L$--operators starting from the Yang's solution of the Yang-Baxter equation
 without referring to the precise meaning of the $L$-operators and to their
 analytical properties.
The result can be thought of as a formal algebra which turns to be genuine
 Hopf algebra when one defines the $L$-operators to be explicit
 generating functions of its elements.

 Let $\overline{R}(u)$ be a rational solution of the quantum Yang--Baxter
 equation:
\beq
\overline R (u)\ =  \  {u-i\h P\over u-i\h}\ ,\label{R-mat}
\eeq
where $P\in {\mbox{End}}\,\CC^2\otimes \CC^2$,
$$P
=\left(\begin{array}{cccc} 1&0&0&0\\ 0&0&1&0\\ 0&1&0&0\\ 0&0&0&1
\end{array}\right)
$$ is a permutation operator and $\h$ is the deformation
parameter. Note that the algebras which we  discuss
below can be defined for arbitrary complex values of the
deformation parameter. But the representation theory, especially the
infinite-dimensional representation theory, depends on the specific
values of this parameter. So we fix $\h\in\RR$ and $\h>0$. Moreover,
this choice is in accordance with applications to massive field theory,
where one should put
 $\h=\pi$ \cite{L,S1}.
%%%%%%%%%%%%%%%%%%%%%%%%%%%%%%%%%%%%%%%%%%%%%%%%%%%%%%%%%

Following the formalism of Faddeev--Reshetekhin--Takhtadjan
\cite{FRT} one can
 use $\overline{R}(u)$ for the construction of the bialgebra whose generators
 are gathered into the matrix elements of the quantum $L$-operator $L(u)$,
\beq
L(u)=\left(\begin{array}{cc}
L_{11}(u)&L_{12}(u)\\ L_{21}(u)&L_{22}(u)
\end{array}\right)
\label{L-op}
\eeq
which satisfy the Yang-Baxter relation
$$\overline{R}(u-v)L_1(u)L_2(v)=L_2(v)L_1(u)\overline{R}(u-v)\ .$$
Here $L_1(u)=L(u)\ot1$ and $L_2(u)=1\ot L(u)$.

More precisely, we would like to consider the Hopf algebras, which are
 quantum doubles
and which have a family of $(2\ell+1)$-dimensional representations
$\pi^{(\ell)}_z$ parametrized by
the parameter $z\in\CC$. Here $\ell$ is the spin of the representations.
For the simplest nontrivial two-dimensional representation we require
 the $L$--operator to be proportional to the $R$-matrix \r{R-mat}:
\beq
\pi^{(1/2)}_z \left(L(u)\right) \sim \overline R(u-z)
\ .\label{2dim}
\eeq

%For the definition of the algebra of $L$-operators we need the
%$R$-matrices $R^\pm(u)$ which differ from the rational
%$R$-matrix \r{R-mat} by the scalar factors: $R^\pm(u)=\rho^\pm(u)
%\overline R(u)$, where

Due to \cite{FRT} such an algebra can be defined via two generating
 matrix-valued functions $L^\pm(u)$ which satisfy the relations
\begin{eqnarray}
\overline R(u_1-u_2)L^+_1(u_1)L^-_2(u_2)&=&
L^-_2(u_2)L^+_1(u_1) \overline R(u_1-u_2)\ ,\nn\\
\overline R(u_1-u_2)L^\pm_1(u_1)L^\pm_2(u_2)&=&
L^\pm_2(u_2)L^\pm_1(u_1) \overline R(u_1-u_2)\ ,
\label{c=0}\\
\qdet\, L^\pm(u)=
L^\pm_{11}(u-i\h)L^\pm_{22}(u)&-&L^\pm_{12}(u-i\h)L^\pm_{21}(u)
=1\ .\label{qdet}
\eea
As usual, the relation \r{qdet} on the q-determinant
factorizes the algebra over
primitive central elements which appear due to degeneracy of the $R$-matrix
at critical points $u=\pm i\h$. This ensures the existence of antipode and
 turn the bialgebra into a Hopf algebra. Two $L$-operators $L^\pm(u)$ generate
 two dual Hopf subalgebras of the quantum double.

It was shown in \cite{KT} that
the relations \r{c=0} \r{qdet} can be interpreted as the defining relations
 for the
 quantum double of the
Yangian associated with $\frak{sl}_2$ \cite{Dr}.

The two-dimensional representation \r{2dim} is
\beq
\pi^{(1/2)}_z \left(L^\pm(u)\right) =\rho^\pm(u-z) \overline R(u-z)
\ ,\label{2dimpm}
\eeq
where the functions $\rho^\pm(u)$ satisfy the equation
$$\rho^\pm(u)\rho^\pm(u-i\h)={u-i\h\over u}\ ,$$
which follows from $q$-determinant condition
 \r{qdet}
and could be chosen, in particular, as follows:
\beq
\rho^\pm (u)=
\left[{\Gamma^2\left(
\frac{1}{2}\mp \frac{u}{2i\h}\right)
\over
\Ga{1\mp\frac{u}{2i\h}} \Ga{\mp\frac{u}{2i\h} } }\right]^{\mp1}
.
\label{rho}
\eeq
Moreover, the solution \r{rho} is unique for certain analyticity conditions
 on the $L$-operators $L^\pm(u)$ which we discuss later. Nevertheless
the formal current algebras which we define further do not depend on the
choice of $\rho^\pm(u)$.

The latter algebra of $L$-operators appears to have only finite-dimensional
 rep\-re\-sen\-ta\-tions.
In order to construct infinite-dimensional representations
 we need to perform the central extension of this formal algebra. It can be
done as follows.
The algebra of $L$-operators \r{c=0}, \r{qdet} admits a family of shifting
 automorphisms
$T_z L^\pm(u)=
L^\pm(u-z)$.
Then the extension of the quantum double by the infinitisemal shift
operator $d$ and the dual to it element $c$ \cite{Drnew,K} leads
to the following commutation relations \cite{RS,FR}:
\begin{eqnarray}
[L(u),c{]}&=&0\ ,\nn\\
\ee^{ad}L^{\pm}(u)&=&L^\pm(u+a)\ee^{ad},\nn\\
R^+(u_1-u_2+ic\tih/2)L^+_1(u_1)L^-_2(u_2)&=&
L^-_2(u_2)L^+_1(u_1) R^+(u_1-u_2-ic\tih/2),\nn\\%\label{RLL-univ-+-}\\
R^\pm(u_1-u_2)L^\pm_1(u_1)L^\pm_2(u_2)&=&
L^\pm_2(u_2)L^\pm_1(u_1) R^\pm(u_1-u_2),
\label{RLL-univ}\\
\qdet\, L^\pm(u)=
L^\pm_{11}(u-i\h)L^\pm_{22}(u)&-&L^\pm_{12}(u-i\h)L^\pm_{21}(u)
=1\ ,\nn
\eea
where
$$R^\pm(u)=\rho^\pm(u)\overline R(u)\ .$$
 We call this algebra the (central extended) algebra of $L$-operators.
 Let us stress that this is still a formal algebra since we did not specify
 the meaning of its generating functions $L^\pm(u)$.
%%%%%%%%%%%%%%%%%%%%%%%%%%%%%%%%%%%%
%%%end of proofs%%
%%%%%%%%%%%%%%%%%%%%%%%%%%%%%

 We can go further to proceed in algebraic manipulations with formal
 algebra of $L$-operators.
 First, it is natural to factorize explicitly the quantum determinat
 and have three generating functions instead of four analogous to classical
 passage from $\frak{gl}_2$ to $\frak{sl}_2$. This can be done by means of
 Gauss decomposition of the $L$-operators.

Let
\beq
L^\pm(u)
=\left(\begin{array}{cc} 1& f^\pm(u)\\0&1\end{array}\right)
\left(\begin{array}{cc}  k^\pm_1(u)&0\\ 0&k^\pm_2(u) \end{array}\right)
\left(\begin{array}{cc} 1&0\\ e^\pm(u)&1\end{array}\right)\ ,
\label{GL-univ}
\eeq
be the Gauss decomposition of $L$-operators.
The condition \r{qdet} implies that the product of the entries
 of diagonal matrix in r.h.s. of \r{GL-univ} is equal to one:
$k^\pm_1(u)k^\pm_2(u+i\h)=1$, therefore they
are invertible and
$$
k^\pm_1(u)=(k^\pm_2(u+i\tih))^{-1}.
$$
Let
$$%\beq
h^\pm(u)=\sk{k^\pm_2\left(u+i\h\right)}^{-1} \sk{k^\pm_2\left(u\right)}^{-1}
\ .
$$%\label{h-new}\eeq
We call the operator valued generating functions
 $e^\pm(u)$, $f^\pm(u)$ and  $h^\pm(u)$
the Gauss coordinates of the $L$-operators.
 We have the following

\begin{prop}\label{3}
 {\sl The Gauss coordinates $e^\pm(u)$, $f^\pm(u)$ and $h^\pm(u)$
satisfy the following commutation relations }
\bea
h^\pm(u)h^\pm(v)&=&h^\pm(v)h^\pm(u)\ ,\nn\\
{[}e^\pm(u),f^\pm(v){]}
&=&i\h {h^\pm(u)-h^\pm(v)\over u-v}\ ,\nn\\
{[}h^\pm(u),e^\pm(v){]}
&=&i\h {\{h^\pm(u),e^\pm(u)-e^\pm(v)\}\over u-v}\ ,\nn\\
{[}h^\pm(u),f^\pm(v){]}
&=&-i\h {\{h^\pm(u),f^\pm(u)-f^\pm(v)\}\over u-v}\ ,\nn\\
{[}e^\pm(u),e^\pm(v){]}
&=&i\h {(e^\pm(u)-e^\pm(v))^2\over u-v}\ ,\nn\\
{[}f^\pm(u),f^\pm(v){]}
&=&-i\h {(f^\pm(u)-f^\pm(v))^2\over u-v}\ ,\nn\\
h^+(u)h^-(v)&=&
{
(u-v-i\tih(1+c/2))
(u-v+i\tih(1+c/2))
\over
(u-v+i\tih(1-c/2))
(u-v-i\tih(1-c/2))
}\
h^-(v)h^+(u)\ ,\nn\\
{[}e^\pm(u),f^\mp(v){]}
&=&
i\h {h^\pm(u)\over u-v\pm ic\h/2}
-
i\h {h^\mp(v)\over u-v\mp ic\h/2}
\ ,\nn\\
{[}h^\pm(u),e^\mp(v){]}
&=&i\h {\{h^\pm(u),e^\pm(u)-e^\mp(v)\}\over u-v\mp ic\h/2}\ ,\nn\\
{[}h^\pm(u),f^\mp(v){]}
&=&-i\h {\{h^\pm(u),f^\pm(u)-f^\mp(v)\}\over u-v\pm ic\h/2}\ ,\nn\\
{[}e^\pm(u),e^\mp(v){]}
&=&i\h {(e^\pm(u)-e^\mp(v))^2\over u-v\mp ic\h/2}\ ,\nn\\
{[}f^\pm(u),f^\mp(v){]}
&=&-i\h {(f^\pm(u)-f^\pm(v))^2\over u-v\pm ic\h/2}\ .\nn
\eea
\end{prop}

\noindent
{\sl Proof.}
The proof is a direct substitution of the Gauss decomposition
of $L$-operators \r{GL-univ} into \r{RLL-univ}. The particular
 case of \r{RLL-univ} at $u_1=u_2-i\h$:
\bea
k^\pm_2(u)e^\pm(u)&=&e^\pm(u-i\h)k^\pm_2(u)\ ,\nn\\
k^\pm_2(u)f^\pm(u-ih)&=&f^\pm(u)k^\pm_2(u)\ ,\nn
\eea
is very useful for this algebraic exercise.
 For the analogous treatment of quantum affine algebras see \cite{DF}.

One can see that the algebra of Gauss coordinates does not depend
 on a choice of the factor $\rho^\pm(u)$ in the definition of $R^\pm(u)$.
 It refers only to the original Yang $R$-matrix.
\bigskip

\noindent
{\bf \razdel Hopf structure of the algebra of $L$-operators.}
The $L$-operator's language is convenient for the description of the
 coalgebraic structure.
The comultiplication map for the algebra \r{RLL-univ}
 of formal $L$-operators is given by the formulas
\bea
\Delta\,c&=&c^{(1)}+c^{(2)}=c\otimes1+1\otimes c\ ,\nn\\
\Delta\,d&=&d\otimes1+1\otimes d\ ,\nn\\
\Delta' L^\pm(u)&=&L(u\pm i\tih c^{(2)}/4)
{\dot\otimes} L(u\mp i\tih c^{(1)}/4)
\label{comul-L-univ}
\eea
or in components
\beq
\Delta L^\pm_{ij}(u)=
\sum_{k=1}^2
L^\pm_{kj}(u\mp i\tih c^{(2)}/4) \otimes L_{ik}(u\pm i\tih c^{(1)}/4)\ .
\label{com-L*cmp}
\eeq
The antipode and counit are:
\bea
S(L^\pm(u))&=&(L^\pm(u))^{-1}\ ,\label{antipode}\\
\epsilon(L^\pm_{ij}(u))&=&\delta_{ij}\ .\label{counit}
\eea

The comultiplications of the Gauss coordinates $e^\pm(u)$, $f^\pm(u)$
and $h^\pm(u)$ reads as follows:
\bea
\Delta e^\pm(u) &=&
e^\pm(u')\ot 1 +
\sum_{p=0} ^{\infty}(-1)^p
       \left(f^\pm(u'-i\tih)\right)^{p} h^\pm(u')\otimes
\left(e^\pm(u'') \right)^{p+1},   \label{com-e-fu}\\
\Delta f^\pm(u)&=&
1\otimes f^\pm(u'') +
\sum_{p=0} ^{\infty} (-1)^p
       \left(f^\pm(u')\right)^{p+1} \otimes
h^\pm(u'')\left(e^\pm(u''-i\tih)
\right)^p,   \label{com-f-fu}\\
\Delta h^\pm(u)&=&\sum_{p=0}^\infty (-1)^p (p+1)
\left(f^\pm(u'-i\tih)\right)^{p} h^\pm(u') \ot
h^\pm(u'')\left(e^\pm(u''-i\tih) \right)^p,
\label{comul-h}
\eea
where $u'=u\mp i\tih c^{(2)}/4$ and $u''=u\pm i\tih c^{(1)}/4$.
The proof of these formulas   can be found in \cite{KLP3}.
\bigskip

\noindent
{\bf \razdel  Current realization of the algebra of $L$-operators.}
To construct the infinite dimensional
rep\-re\-sen\-ta\-tion theory of the central extended $L$-operator algebra
the Gauss coordinates are
inconvenient and it is natural to introduce the total currents.
Let
\bea
e(u)&=&
e^+\sk{u+\fract{ic\tih}{4}}-e^-\sk{u-\fract{ic\tih}{4}}\ ,\nn\\
f(u)&=&
f^+\sk{u-\fract{ic\tih}{4}}-f^-\sk{u+\fract{ic\tih}{4}}\ ,\label{total}
\eea
be the combination of the Gauss coordinates of the $L$-operators
which we call the total currents.  One can verify that the
commutation relations between the total currents and $h^\pm(u)$ close:
\bea
[d,e(\u){]}&=&\frac{{d}}{{d}u}e(\u), \hsp
 [d,f(\u)]=\frac{{d}}{{d}u}f(\u), \nn \\
h^\pm(u)e(v)&=&
{ (u-v-i\tih(1\pm c/4))
\over
(u-v+i\tih(1\mp c/4))
}\
e(v)h^\pm(u)\ , \nn\\
h^\pm(u)f(v)&=&
{ (u-v+i\tih(1\pm c/4))
\over
(u-v-i\tih(1\mp c/4))
}\
f(v) h^\pm(u)\ , \nn\\
e(u)e(v)&=&
{ (u-v-i\tih)
\over (u-v+i\tih)
}\  e(v)e(u)\ ,\nn\\
f(u)f(v)&=&
{ (u-v+i\tih)
\over (u-v-i\tih)
}\ f(v)f(u)\ ,\nn
\eea
\beq%$$%\bea
{[}e(u),f(v){]}= -i\h\left[
\delta\left(u-v+\fract{ic\tih}{2}\right)
h^+\left(u+\fract{ic\tih}{4}\right)-
\delta\left(u-v-\fract{ic\tih}{2}\right)
h^-\left(v+\fract{ic\tih}{4}\right)
\right]  ,       \label{CC1}
\eeq%$$%\eea
where the $\delta$-function is defined by the formal equality
$$ \delta(u-v)=\frac{1}{u-v}-\frac{1}{u-v}$$
and satisfies the relation
$$g(u)\delta(u-v)=g(v)\delta(u-v)\ .$$
In the next section we develop a technique of extracting the
 correctly define (Hopf) algebra from the above current algebra
provided certain analytical conditions are supposed.

\setcounter{equation}{0}
\setcounter{parag}{0}
\section{Factorization and the Riemann problems}

{\bf \razdel}
Until now we did not specify the analytical properties of $L$-operators.
Let us fix them in the following manner:
\bea
&L^+(u)\quad \hbox{is analytical in some neighborhood of}\ u=\infty\ , \nn\\
&L^-(u)\quad \hbox{is analytical in some neighborhood of}\ u=0\ . \nn
\eea
To simplify the consideration  we set $c=0$ first and restore the
central element only in the final formulas. We would like to invert
the relation
\beq
e(u)=e^+(u)-e^-(u)\ ,\label{Soh1}
\eeq
namely to express the generating functions $e^\pm(u)$ through the total
current $e(u)$  by means of some integral transforms. More precisely,
 let us suppose that

\begin{itemize}
\item[{$(i)$}]the function $e(u)$ is analytical function which has only
 isolated singularities on compactified complex plane;
\item[{$(ii)$}]the functions $e^+(u)$ is analytical at infinity and
 the function $e^-(u)$ is analytical at zero.
\end{itemize}

We call the solution of \r{Soh1} with the conditions
$(i)$ and $(ii)$ the classical Riemann problem for a circle.
Let us first
fix a closed counterclockwise oriented contour $\Gamma$  around $0$
 (for instance, $\Gamma=\{|u|=1\}$). Then the
 Cauchy type integrals $\tilde
e^\pm(u) = \oint_\Gamma {dv\over 2\pi i} {e(v)\over
u-v}$, where $u$ is outside or inside $\Gamma$ give two functions $\tilde
e^\pm(u)$
 analytical outside and inside of $\Gamma$ such that the relation \r{Soh1}
 is valid for the points $u\in\Gamma$. The limiting values of these functions
 when $u$ tends to a contour are given by Sokhotsky-Plemely relations:
\beq
\tilde e^\pm(u)={1\over 2}\left[\mbox{V.P.}
\oint_\Gamma{dv\over \pi i}{e(v)\over u-v} \pm e(u)\right],\quad
u\in\Gamma\ . \label{Soh3}
\eeq

Suppose also that the function $e(u)$ is analytical in the points  of
$\Gamma$ so we have  the identity
\beq
e(u)=\pm \int_{C^\pm_{u,\epsilon}} {dv\over \pi i} {e(v)\over u-v}\ ,
\label{id-Cauchy}
\eeq
where  $C^\pm_{u,\epsilon}$ are small semicircles around the point $u$
of the radius $\epsilon$ drown in the different directions
as shown on the Fig.~1.

\unitlength 1mm
\linethickness{0.4pt}
\begin{picture}(120.00,33.00)
\bezier{76}(70.00,21.00)(71.00,30.67)(80.00,31.00)
\bezier{76}(80.00,31.00)(89.00,31.33)(90.00,21.00)
\put(90.00,21.00){\vector(0,-1){0.2}}
\bezier{76}(70.00,19.00)(70.33,9.33)(80.00,9.00)
\bezier{76}(80.00,9.00)(89.67,10.00)(90.00,19.00)
\put(90.00,19.00){\vector(0,1){0.2}}
\put(80.00,20.00){\makebox(0,0)[cc]{$\cdot$}}
%\put(86.00,29.00){\vector(2,3){0.2}}
\put(82.00,21.33){\makebox(0,0)[cc]{$u$}}
\put(73.00,33.00){\makebox(0,0)[cc]{$C^+_{u,\epsilon}$}}
\put(72.67,6.00){\makebox(0,0)[cc]{$C^-_{u,\epsilon}$}}
\put(120.00,5.00){\makebox(0,0)[cc]{Fig.~1.}}
\end{picture}

\noindent
Using \r{id-Cauchy} we can
obtain the solution of \r{Soh1} in the domain of analyticity of $e(u)$
nearby the contour $\Gamma$:
\beq
e^\pm(u)=
\oint_{\Gamma_\pm}{dv\over 2\pi i}{e(v)\over u-v}=
\oint_{|v|\lessgtr |u|}{dv\over 2\pi i}{e(v)\over u-v}\ ,
\label{gamma}
\eeq
where the contours $\Gamma_\pm=\Gamma\oplus C^\pm_{u,\epsilon}$
are  shown on the Fig.~2. The distinction between the generating
functions $\tilde e^\pm(u)$ and $e^\pm(u)$ is that the former are
related to the Riemann problem with fixed contour, while the latter
with the problem where contours $\Gamma_\pm$ are not strictly fixed.
\medskip

\unitlength 1mm
\linethickness{0.4pt}
\begin{picture}(120.00,45.00)
\bezier{84}(90.00,25.00)(90.33,35.33)(80.00,35.00)
\bezier{76}(80.00,35.00)(71.00,35.33)(70.00,25.00)
\bezier{80}(70.00,25.00)(69.67,15.33)(80.00,15.00)
\bezier{76}(80.00,15.00)(90.33,16.00)(90.00,25.00)
\put(90.00,26.00){\vector(0,1){0.2}}
\bezier{108}(81.00,10.00)(94.00,11.00)(95.00,25.00)
\bezier{116}(95.00,25.00)(94.67,39.33)(80.00,40.00)
\bezier{112}(80.00,40.00)(66.00,39.33)(65.33,25.00)
\bezier{112}(65.33,25.00)(66.00,10.33)(79.00,10.00)
\put(79.00,10.00){\vector(1,0){0.2}}
\put(81.00,10.00){\vector(-1,0){0.2}}
\put(80.00,10.00){\makebox(0,0)[cc]{$\cdot$}}
\put(80.00,25.00){\makebox(0,0)[cc]{$\cdot$}}
\put(80.00,22.00){\makebox(0,0)[cc]{$0$}}
\put(65.33,24.00){\vector(0,-1){0.2}}
\bezier{160}(100.00,25.00)(100.00,45.00)(80.00,45.00)
\bezier{156}(80.00,45.00)(60.67,44.33)(60.00,25.00)
\bezier{156}(60.00,25.00)(60.67,5.00)(80.00,5.00)
\bezier{160}(80.00,5.00)(100.00,5.00)(100.00,25.00)
\put(100.00,26.00){\vector(0,1){0.2}}
\put(63.67,25.00){\makebox(0,0)[cc]{$\Gamma$}}
\put(87.33,25.00){\makebox(0,0)[cc]{$\Gamma_+$}}
\put(105.00,25.00){\makebox(0,0)[cc]{$\Gamma_-$}}
\put(120.00,5.00){\makebox(0,0)[cc]{Fig.~2.}}
\put(80.00,8.20){\makebox(0,0)[cc]{$u$}}
\end{picture}

\noindent
Now we forget the initially fixed contour $\Gamma$ and use the
 integral transforms \r{gamma} as the solution of \r{Soh1}
 satisfying the condition $(ii)$, which is valid due to $(i)$.
The contours $\Gamma_+$ ($\Gamma_-$) should be chosen close to the
point $u$ which means that they are located in the regipns
$|u|-\epsilon <|v|<|u|$ ($|u|<|v|<|u|+\epsilon$)
of the analyticity of the functions $e(v)$.

Restoring the dependence on the central element in \r{total} we obtain
 analogously:
\bea
e^+(u)&=&\oint_{|v|<|u-ic\h/4|} {dv\over2\pi i}\
 { e(v)\over
(u-v-ic\tih/4)}\ ,\nn\\
f^+(u)&=&\oint_{|v|<|u+ic\h/4|} {dv\over2\pi i}\
 { f(v)\over
(u-v+ic\tih/4)}\ ,\nn\\
e^-(u)&=&\oint_{|v|>|u+ic\h/4|} {dv\over2\pi i}\
 { e(v)\over
(u-v+ic\tih/4)}\ ,\nn\\
f^-(u)&=&\oint_{|v|>|u-ic\h/4|} {dv\over2\pi i}\
 { f(v)\over
(u-v-ic\tih/4)}\ ,\label{int-d}
\eea
where all contours in the integrals are
counterclockwise circles around the point $v=0$
such that the corresponding points $u\pm i\h c/4$ are either outside
of the contours for the generating functions $e^+(u)$, $f^+(u)$ or inside
the contours for the generating functions $e^-(u)$, $f^-(u)$.

The integral transformations \r{int-d}
are the solutions of the Riemann problems
for the factorization of the entire function into sum of the functions
analytical in the neighborhood of the given point ($u=0$)
and $u=\infty$ respectively.

The $\delta$-function in the commutation relation of the total
currents $e(u)$ and $f(v)$ should be also presented in terms of the same
 Riemann problem:
\beq
\delta(u-v)=\left.{1\over u-v}\right|_{|u|>|v|} -
\left.{1\over u-v}\right|_{|u|<|v|}=
\sum_{n+m=-1}u^nv^m
\ .
\label{delta}
\eeq

The relations \r{int-d} dictate the precise sense of the
operator valued functions $e(u)$, $f(u)$, $e^\pm(u)$ and $f^\pm(u)$
as the generating series of the elements of the algebra $\DY$ (central
extended Yangian double).
%%%%%%%%%%%%%%%%%%%%%%%
If we decompose the currents $e^\pm(u)$ into Taylor series in the
 points of their regularity ($\infty$ and $0$):
$$%\beq\label{2.0}
e^{\pm}(\u)=\mp
i\h\sum_{k\geq0\atop k<0} e_k (\u\mp ic\h/4)^{-k-1}\ , \quad f^{\pm}(\u)=\mp
i\h\sum_{k\geq0\atop k<0} f_k (\u\pm ic\h/4)^{-k-1}\ ,  $$
then, due to \r{Soh1}, we have the presentation
\beq
e(u)=-i\h \sum_{n\in\ZZ} e_n u^{-n-1}\ ,\quad
f(u)=-i\h \sum_{n\in\ZZ} f_n u^{-n-1}\ .
\label{total-d}
\eeq
 We see that the decompositions \r{total}
 modulo technical shifts are just the decompositions of
a formal power series into the parts with positive and negative powers
 and the coefficients of the series are given as
\beq
\label{456}
e_n=\pm{1\over2\pi\h}\oint_{C_\pm}v^ne^\pm(v\mp ic\h/4)dv,
\qquad  f_n=\pm{1\over2\pi\h}\oint_{C_\pm}v^nf^\pm(v\pm ic\h/4)dv\ ,
\eeq
or, equivalently (this follows from the definition \r{gamma})
\beq
\label{457}
e_n={1\over2\pi\h}\oint_{C_\pm}v^ne(v)dv,
\qquad  f_n={1\over2\pi\h}\oint_{C_\pm}v^nf(v)dv\ ,
\eeq
 where $C_+$ is a contour surrounding the infinity
(clockwise, like $\Gamma_+$) and is taken
for $n\geq0$ while
$C_-$ is a contour
 surrounding zero (counterclockwise, like $\Gamma_-$)
and is taken for $n<0$.
Note that in the following we use \r{456} as basic definitions since for
the functions $e^\pm(u)$ and $f^\pm(u)$ of the distribution
type, the integrals in \r{456} make precise sense contrary to
the integrals in \r{457}.
It is natural also to assume that the currents $h^\pm(u)$ satisfy
 the same analyticity conditions as $e^\pm(u), f^\pm(u)$ and put
 $$%\beq
h^\pm(u)=1\mp i\h\sum_{k\geq0\atop k<0}h_{k}\u^{-k-1},
$$%\label{yancurr}\eeq
where
\beq
h_n=\pm{1\over2\pi\h}\int_{C_\pm}(h^\pm(u)-1) u^n du
\label{458}
\eeq
 with the same rule for the signs.
The unity in the above formula is introduced because of the group-like
nature of the Cartan generating series.
\bigskip

\noindent
{\bf \razdel}
There is another choice of the analytical data. Let us fix the following
analytical behavior of the operators $L^\pm(u)$:
\bea
&L^+(u)\quad \hbox{is analytical in $\CC$ for}\quad \Im u<-A\ , \nn\\
&L^-(u)\quad \hbox{is analytical in $\CC$ for}\quad \Im u>A\ , \nn
\eea
for some positive $A$.
In an analogous manner we obtain that for a meromorphic function $e(u)$
decreasing when $ \Re u \to \pm\infty$
the Cauchy type integrals
$$ \tilde e^\pm(u)=
\int_{-\infty}^\infty{dv\over 2\pi i}{e(v)\over u-v}
$$
 are analytical in closed half-planes
$\Im  u \leq 0$ and  $\Im u \geq 0$;
they satisfy the relation \r{Soh1} for real
$u$ and can be used for the following presentations of the functions
$e^\pm(u)$ satisfying \r{Soh1} for all $u$ and
analytical in half-planes $\mbox{Im}\ u <-A$ and
 $\mbox{Im}\ u >A$ for some positive $A$:
\beq
e^\pm(u)=
 \int_{\Im v\ \gtrless\ \Im u}
{dv\over 2\pi i}{e(v)\over u-v}\ .
\label{311}
\eeq
The contours of the integrations in \r{311} are close to the point
$u$ which means that there are no singularities of the function
$e(z)$ in the strips $\Im u \lessgtr \Im z \lessgtr \Im v$.

Restoring the central elements we obtain:
\bea
e^+(u)&=&\int_{\Im (u-v)<c\h/4}
{dv\over2\pi i}\
 { e(v)\over
(u-v-ic\tih/4)},\nn\\
f^+(u)&=&\int_{\Im (u-v)<-c\h/4}
{dv\over2\pi i}\
 { f(v)\over
(u-v+ic\tih/4)}\ ,\nn\\
e^-(u)&=&\int_{\Im (u-v)>-c\h/4}
{dv\over2\pi i}\
 { e(v)\over
(u-v+ic\tih/4)},\nn\\
f^-(u)&=&\int_{\Im (u-v)>c\h/4}
{dv\over2\pi i}\
 { f(v)\over
(u-v-ic\tih/4)}\ .\label{int-f}
\eea
The corresponding integral transforms are related to the Riemann
problem of fac\-to\-riz\-ing the decreasing meromorphic function into a sum of
functions analytical in certain
 upper and lower half-planes of the complex plane.

The $\delta$-function in the commutation relation of the total
currents $e(u)$ and $f(v)$ should be also the solution
 of the same Riemann problem and
is given by the
integral:  \beq \delta(u-v)=\lim{\epsilon\to0} \left[ {1\over
u-v-i\epsilon}-{1\over u-v+i\epsilon}\right]= i \intt d\la\ \ee^{-i\la
(u-v)}\ .
\label{delta-con}
\eeq

Analogously to the case of Riemann problem for a circle,
the integral relations \r{int-f} dictate the precise sense of the
operator valued functions  $e(u)$, $f(u)$, $e^\pm(u)$ and $f^\pm(u)$
as the generating integrals of the elements of the algebra
$\Aya$ (degeneration of $\Ael$, $\eta\to0$, \cite{KLP3}). Namely,
the analytical in half-planes $\Im\,u<-A$ and $\Im\,u>A$ functions
$e_+(u)$, $f_+(u)$ and $e_-(u)$, $f_-(u)$ can be presented
via Laplace integrals:
\beq
e^\pm(u)=\tih \int_{0}^{\pm\infty} d\la\ \ee^{-i\la u}\
\he_\la \ee^{-c\tih|\la|/4}\ ,\quad  %\label{Lapl-e}\\
f^\pm(u)\ =\ \tih \int_{0}^{\pm\infty} d\la\ \ee^{-i\la u}\
\hf_\la \ee^{c\tih|\la|/4}\ . \label{Lapl-f}
%\ka^\pm(u)&=&\int_{0}^{\pm\infty} d\la\ \ee^{i\la u}\
%\kk_\la \ ,\label{Lapl-kk}\\
\eeq
Then, due to \r{total}, we have the presentation of the
total currents
\beq
e(u)=\h\int_{-\infty}^{\infty}d\la\ \ee^{-i\la u} \he_\la\ ,\quad
f(u)=\h\int_{-\infty}^{\infty}d\la\ \ee^{-i\la u} \hf_\la\
\label{total-f}
\eeq
and the relation \r{total} is the decomposition of formal Fourier integral
into the sum of two Laplace transforms.

Due to the definition and the invertion formulas for the Laplace
transforms, the Fourier modes $\he_\la$ and $\hf_\la$ are the following
integrals:
\bea
\he_\la&=&\pm{\ee^{c\h|\la|/4}\over2\pi \h}
\int_{C_\pm}du\,\ee^{i\la u}e^\pm(u),\quad
\he_0={1\over 2\pi \h}\left(\int_{C_+}du\, e^+(u)-  \int_{C_-}du\, e^-(u)
\right),\label{676}\\
\hf_\la&=&\pm{\ee^{-c\h|\la|/4}\over2\pi \h}\int_{C_\pm}du\,\ee^{i\la
u}f^\pm(u),\quad \hf_0={1\over 2\pi \h}\left(\int_{C_+}du\, f^+(u)-
\int_{C_-}du\, f^-(u) \right),\label{678}
\eea
where $C_+$ is the line parallel
to the real axis and belonging to the half-plane $\Im\,u<-A$ for $\la>0$ and
$C_-$ is also the line parallel to the real axis but in another
half-plane $\Im\,u>A$ for $\la<0$. All the integrals are principal
 value integrals. Again we can use total currents
$e(u)$ and $f(u)$ in  \r{678}:
\bea
\he_\la&=&{1\over2\pi \h}
\int_{C_\pm}du\,\ee^{i\la u}e(u),\quad
\he_0={1\over 2\pi \h}\left(\int_{C_+}du\, e(u)+  \int_{C_-}du\, e(u)
\right),\nn\\
\hf_\la&=&{1\over2\pi \h}\int_{C_\pm}du\,\ee^{i\la
u}f(u),\quad \hf_0={1\over 2\pi \h}\left(\int_{C_+}du\, f(u)+
\int_{C_-}du\, f(u) \right).\nn
\eea
It is natural to suppose also that
\bea
h^\pm(u)= %\exp\left(2\tih\ka^\pm(u) \right)=
1+\tih \int_{0}^{\pm\infty} d\la\ \ee^{-i\la u}\
\hhh_\la \ee^{-c\h|\la|/4}
\ ,\label{Lapl-h}
\eea
where $\hhh_\la$ are given by principal value integrals
\bea
\hhh_\la&=&\pm{\ee^{c\h|\la|/4}\over2\pi \h}
\int_{C_+}du\,\ee^{i\la u}(h^+(u)-1),\quad
\mbox{for}\quad \la\not=0\ ,\label{677}\\
\hhh_0&=&{1\over2\pi \h}\left(
\int_{C_+}du\,(h^+(u)-1)-
\int_{C_-}du\,(h^-(u)-1)\right).
\label{679}
\eea

We see that two different Riemann problems correspond to two different
algebras. We show in the next sections that they are the central
extended Yangian double and the rational degeneration of $\Ael$.
 Let us note once more that these Hopf algebras can be
completely defined by the formal algebra of
total currents $e(u)$, $f(u)$, $h^\pm(u)$
and by the type of the Riemann problem. We will return to this point
 in the last section.

\setcounter{equation}{0}
\setcounter{parag}{0}
\section{Algebras $\DYsect$ and $\Ayasect$}

{\bf \razdel The algebra $\DY$.}
Let us consider the above algebra specialized by the Riemann
  problem \r{int-d} for a circle in more details.
For the readers convenience we identify $-i\h=\yh$.
 As we have seen before, this algebra
 (we denote it for a short time by letter $D$) can be defined as an algebra
 generated by the elements $c, d, e_n, f_n, h_n, n\in \ZZ$ gathered
 into generating functions
 $$
e(u)=\yh\sum_{k\in\ZZ}e_ku^{-k-1},\qquad
f(u)=\yh\sum_{k\in\ZZ}f_ku^{-k-1},$$
\beq
h^\pm(u)=1\pm\yh\sum_{k\geq0\atop k<0}h_{k}\u^{-k-1}
\label{yancurr}
\eeq
which satisfy the relations
%%%%%%%%%%%
\bea
[d,e(\u){]}&=&\frac{\mbox{d}}{\mbox{d}u}e(\u), \hsp
 [d,f(\u)]=\frac{\mbox{d}}{\mbox{d}u}f(\u), \nn \\
h^+(u)h^-(v)&=&
{
(u-v+\yh(1+c/2))
(u-v-\yh(1+c/2))
\over
(u-v-\yh(1-c/2))
(u-v+\yh(1-c/2))
}\
h^-(v)h^+(u)\ ,\nn\\
h^\pm(u)e(v)&=&
{ (u-v+\yh(1\pm c/4))
\over
(u-v-\yh(1\mp c/4))
}\
e(v)h^\pm(u)\ , \nn\\
h^\pm(u)f(v)&=&
{ (u-v-\yh(1\pm c/4))
\over
(u-v+\yh(1\mp c/4))
}\
f(v) h^\pm(u)\ , \nn\\
e(u)e(v)&=&
{ (u-v+\yh)
\over (u-v-\yh)
}\  e(v)e(u)\ ,\nn\\
f(u)f(v)&=&
{ (u-v-\yh)
\over (u-v+\yh)
}\ f(v)f(u)\ ,\nn
\eea
\beq
{[}e(u),f(v){]}= \yh\left[
\delta\left(u-v-\fract{c\yh}{2}\right)
h^+\left(u-\fract{c\yh}{4}\right)-
\delta\left(u-v+\fract{c\yh}{2}\right)
h^-\left(v-\fract{c\yh}{4}\right)
\right]  ,            \label{CCeta}
\eeq
%%%%%%%%%%%%%%
%Here we denote $\yh =-i\h$ for the convenience of reading.
The total currents $e(u)$ and  $f(u)$ decompose into the sums
$$e(u)=e^+(u-c\yh/{4})-e^-(u+c\yh/{4})\ ,\quad
 f(u)=f^+(u+c\yh/{4})-f^-(u-c\yh/{4}),$$
where
$$%\beq\label{2.0}
e^{\pm}(\u)=\pm\yh\sum_{k\geq0\atop k<0} e_k (\u\mp c\yh/4)^{-k-1},
\quad
f^{\pm}(\u)=\pm\yh\sum_{k\geq0\atop k<0} f_k (\u\pm c\yh/4)^{-k-1},
$$
due to the  relations \r{int-d} and
the generators of the algebra $e_k, f_k, h_k$ are given by
 the inversion formulas \r{456}, \r{458}.
 The coalgebraic structure is given by the relations
 \r{com-e-fu}--\r{comul-h} %($-i\h=\yh$)
 for the generating functions
$h^\pm(u)$ and $e^\pm(u)$, $f^\pm(u)$.

Let us translate these data into the language of generators
$e_k ,f_k ,h_k$.
Due to the definition of these generators, the translation should be
done in two steps. In the first step, following the definitions
\r{int-d} we get from the relations \r{CCeta} the relations for
the currents $e^\pm(u)$,  $f^\pm(u)$,  $h^\pm(u)$.
These relations will be precise those given by Proposition
\ref{3}. Then we put the spectral parameters into the domains of the
analyticity of these currents and use the definitions
\r{456} and \r{458} of the generators $e_n$, $f_n$ and $h_n$.
An example of such a calculation will be given for a Riemann
problem on a line. Note that for a Riemann problem for a circle the result
coincides with formal substitution of power series
\r{yancurr}, \r{delta} into relations \r{CCeta}.

 At $c=0$ the full set of
the commutation relations is:
\beq
{[}d,e_n{]}=-ne_{n-1}\ ,\quad
{[}d,f_n{]}\ =\ -nf_{n-1}\ ,\quad
{[}d,h_n{]}\ =\ -nh_{n-1}\ ,\label{comm1}
\eeq
\beq
{[}h_0,e_n{]}=2e_{n}\ ,\quad [h_0,f_n{]}\ =\ -2f_{n}\ ,\label{comm2}
\eeq
\beq
{[}h_k,h_n{]}=0\ ,\quad
{[}e_k,f_n{]}=h_{k+n}\ ,\label{comm3}
\eeq
\beq
{[}h_{k+1},e_n{]}-{[}h_k,e_{n+1}{]}=
\yh\{h_k,e_n\}\ ,\quad
{[}h_{k+1},f_n{]}-{[}h_k,f_{n+1}{]}=
-\yh\{h_k,f_n\}\ ,\label{comm4}
\eeq
\beq
{[}e_{k+1},e_{n}{]} - [e_{k},e_{n+1}{]} =
\yh\{ e_{k},e_{n}\}\ ,\quad
[f_{k+1},f_{n}{]} -  [f_{k},f_{n+1}{]} =
-\yh\{f_{k},f_{n}\}\ .\label{commut-dis}
\eeq
When $c\neq0$ the commutation relations of the type $[h_k,e_n{]}$,
$[h_k,f_n{]}$, $[h_k,h_n{]}$ and $[e_k,f_n{]}$ become more
 complicated while the rest are unchanged.
\bea
{[}h_{k+1},e_n{]}&-&{[}h_k,e_{n+1}{]}\ +\ c\yh\vartheta(k)[h_k,e_n]/4\ =\
\yh\{h_k,e_n\}\ ,\label{comm5}\\
{[}h_{k+1},f_n{]}&-&{[}h_k,f_{n+1}{]}\ +\ c\yh\vartheta(k)[h_k,f_n]/4\ =\
-\yh\{h_k,f_n\}\ ,\label{comm6}\\
{[}h_0,h_k{]}&=&0,\quad k\in\ZZ\nn\\
{[}h^+_{k+2},h^-_{-n-1}{]}&-&
2{[}h^+_{k+1},h^-_{-n}{]}
\ +\ {[}h^+_{k+2},h^-_{-n-1}{]}=\nn\\
&=&\eta^2(1+c^2/4){[}h^+_{k},h^-_{-n-1}{]}
-2c\eta^2{\{}h^+_{k},h^-_{-n-1}{\}},\quad k,n\geq 0\ , \label{comm7}\\
{[}e_n,f_p{]}&=&\sum_{k=0}^{n+p}
h_{n+p-k}(-c\yh/4)^{k}B^k_{n,p}\ ,  \label{commut3}\\
{[}e_{-n-1},f_{-p-1}{]}&=&\sum_{k\geq0}
h_{-n-p-k-2}(-c\yh/4)^{k}D^k_{n,p}\ ,  \label{commut2}\\
{[}e_{n},f_{-p-1}{]}&=&\yh^{-1}\sum_{k\geq0}
{[}(-1)^kh^+_{n-p-1-k}-h^-_{n-p-1-k}{]}
(-c\yh/4)^{k}A^k_{n,p}\ ,  \label{commut1}\\
{[}e_{-p-1},f_{n}{]}&=&\yh^{-1}\sum_{k\geq0}
{[}h^+_{n-p-1-k}-(-1)^{k}h^-_{n-p-1-k}{]}
(-c\yh/4)^{k}A^k_{n,p}\ ,  \label{commut}
\eea
where in last four relations $n,p\geq0$,
$$
A^k_{n,p}=\sum_{k'=0}^k C^{k'}_nC_{k+p-k'}^p,\quad
B^k_{n,p}=\sum_{k'=0}^k (-1)^{k'} C^{k'}_nC^{k-k'}_p,\quad
D^k_{n,p}=\sum_{k'=0}^k (-1)^{k'} C^{n}_{n+k'}C^{p}_{p+k-k'},
$$
$C_k^{k'}$ are binomial coefficients,
$$
\vartheta(k)=\left\{ \begin{array}{rl}
1,&\quad k\geq 0\ ,\\ -1,&\quad k< 0\ ,
\end{array}\right.
$$
and in order to write dow the commutation relations between
Cartan generators $h_k$  and the commutation relations
\r{commut1} and \r{commut} we introduced the short notations:
$h^+_{-1}=1$, $h^+_{k}=\yh h_k$, for $k>-1$ and $h^+_k=0$ for
$k<-1$. Analogously, $h^-_{-1}=1-\yh h_{-1}$, $h^-_k=-\yh h_k$
for $k<-1$ and $h^-_k=0$ for $k>-1$.

One can see that the commutation relations \r{comm5}, \r{commut}
coincide with the com\-mu\-ta\-tion relations for the generators of the
central extended Yangian Double $\DY$ (as well as comultiplication rules).
 More precisely, the algebra $D$ defined by the relations \r{commut}
 admit the filtration
\beq
\ldots\subset D_{-n}\subset \ldots\subset D_{-1} \subset D_{0}
\subset D_{1}\ldots\subset D_{n}\ldots \subset D
\label{2.7}
\eeq
defined by conditions $\mbox{deg}\, e_{k}=\mbox{deg}\, f_{k}=
\mbox{deg}\, h_{k}=k$; $\deg \left\{ x\in D_m\right\}\leq m$.
 Then the formal completion of $D$ over this filtration is a Hopf algebra
 which coincides with $\DY$.
\bigskip

\noindent
{\bf \razdel The algebra $\Aya$.}
Let us turn our attention to the Riemann problem for half-planes
 \r{int-d}.
 As we have seen in the previous section, the corresponding algebra
${\cal A}$
 is generated by the elements $\he_\la$, $\hf_\la$, $\hhh_\la$,
 $\la \in\RR$, $c$ and
 $\tilde d$
 gathered into generating integrals
 $$e(u)=\h\int_{-\infty}^{+\infty}d\la\ e^{-i\la u}\he_\la\ ,\qquad
 f(u)=\h\int_{-\infty}^{+\infty}d\la\ e^{-i\la u}\hf_\la\ ,$$
 $$h^\pm(u)=
1+\tih \int_{0}^{\pm\infty} d\la\ \ee^{-i\la u}\
\hhh_\la \ee^{-c\h|\la|/4}\ .$$
 They satisfy the commutation relations \r{CC1} and comultiplication
 rules \r{com-e-fu}--\r{comul-h} for the generating integrals
$h^\pm(u)$ and $e^\pm(u)$, $f^\pm(u)$,  which, due to \r{int-f},
have the form
\bea
e^\pm(u)&=&\tih \int_{0}^{\pm\infty} d\la\ \ee^{-i\la u}\
\he_\la \ee^{-c\tih|\la|/4}\ , \quad%\label{Lapl-e}\\
f^\pm(u)\ =\ \tih \int_{0}^{\pm\infty} d\la\ \ee^{-i\la u}\
\hf_\la \ee^{c\tih|\la|/4}\ , \label{Lapl-e+f}
%\ka^\pm(u)&=&\int_{0}^{\pm\infty} d\la\ \ee^{i\la u}\
%\kk_\la \ .\label{Lapl-kk}\\
\eea
 The elements $\he_\la$, $\hf_\la$, $\hhh_\la$ are given by the inversion
 formulas \r{676}, \r{678}, \r{677}, \r{679}.
%%%%%%%%%%%%%%%%%%%%%%%%%%%%%%
For the description of the algebraical structure of the algebra
$\Aya$ it is more convenient to use  the generators $\kk$
 being the Fourier modes of the logarithm of the currents $h^\pm(u)$:
$$\kappa^\pm(u)=\log h^\pm(u)=\h\int^{\pm\infty}_0d\la\
e^{-i\la u}\kk_\la\ .$$
The generators $\kk_\la$ and $\hhh_\la$ are connected as follows:
\beq
\hhh_\la=\ee^{c\h|\la|/4}
\left(\kk_\la+\sum_{n\geq2}{(\tih)^{n-1}\over n!}
\int_0^\la d\la_1\ \cdots\int_0^{\la_{n-2}}
d\la_{n-1}
\ \kk_{\la_1}\ldots \kk_{\la_{n-1}}
\kk_{\la-\sum_{i=1}^{n-1}\la_i}
\right).
\label{hhh-kk}
\eeq
In particular  $\hhh_0=\kk_0$. The inverse to \r{hhh-kk} relation
is
$$
\kk_\la=\ee^{-c\h|\la|/4}
\left(\hhh_\la+\sum_{n\geq2}{(-\tih)^{n-1}\over n}
\int_0^\la d\la_1\ \cdots\int_0^{\la_{n-2}} d\la_{n-1}
\ \hhh_{\la_1}\ldots \hhh_{\la_{n-1}}
\hhh_{\la-\sum_{i=1}^{n-1}\la_i}
\right).
$$

 The commutation relations given in Proposition \ref{3}
and written in terms of the  generators
 $\he_\la$, $\hf_\la$, $\kk_\la$,
 $\la \in\RR$, $c$ and
 $\tilde d$  take the form:
\bea
[c,\mbox{everything}{]}&=&0,\nn\\
{[}\tilde{d},\he_{\la}{]}\ =\ -\la\he_{\la},\quad
{[}\tilde{d},\hf_{\la}{]}&=& -\la\hf_{\la},\quad
{[}\tilde{d},\kk_{\la}{]}\ =\ -\la\kk_{\la}\ ,\label{d}
\eea
\bea
{[}\kk_\la,\he_\mu{]}&=&{2\sh\,\tih\la\over\tih\la}\ \ee^{-\h c|\la|/4}
\he_{\la+\mu}\ ,\qquad%\label{kk-he}
{[}\kk_\la,\hf_\mu{]}=-{2\sh\,\tih\la\over\tih\la}\ \ee^{\h c|\la|/4}
\hf_{\la+\mu}\ ,\label{kk-hf}\\
{[}\he_\la,\he_\mu{]}&=&\tih \int_{\mu}^{\la}d\tau\
\left[\theta(\tau-\mu)-\theta(\tau-\la)\right]\
\he_{\tau}\he_{\la+\mu-\tau}\ ,\label{he-he}\\
{[}\hf_\la,\hf_\mu{]}&=&-\tih \int_{\mu}^{\la}d\tau\
\left[\theta(\tau-\mu)-\theta(\tau-\la)\right]\
\hf_{\tau}\hf_{\la+\mu-\tau}\ ,\label{hf-hf}\\
{[}\kk_\la,\kk_\mu{]}&=&{4\over\tih^2\la}\ \sh(\tih\la)\,\sh(\tih\la c/2)\,
\delta(\la+\mu)\ ,\label{kk-kk}\\
{[}\he_\la,\hf_\mu{]}&=&{2\tih^{-1}}\
\sh\left({\la\tih c/2}\right)\delta(\la+\mu) +\nn\\
&+&\left[\ee^{(\la-\mu)\tih c/4}\theta(\la+\mu)+
\ee^{(\mu-\la)\tih c/4}\theta(-\la-\mu) \right]
\hhh_{\la+\mu}\ ,\label{he-hf}
\eea
where the step-function $\theta(\la)$ is defined as
$$
\theta(\la)=\left\{\begin{array}{r} 1\quad\mbox{for}\quad \la >0\ ,\\
1/2\quad\mbox{for}\quad \la =0\ ,\\
0\quad\mbox{for}\quad \la<0\ .\end{array}\right.
$$
Note the similarity of the relations \r{kk-hf}, \r{kk-kk} and
\r{he-hf} with some of the relations for the quantum affine
algebra $U_q(\widehat{sl}_2)$ in new realization \cite{Drnew}.

Let us obtain the formula \r{he-he} from the commutation relation
in terms of the total current
$$[e(u),e(v)]=-i\h{\{e(u),e(v)\}\over u-v}\ .$$
The first step is to obtain the commutation relations between
the currents $e^\pm(u)$ using the Riemann problem \r{int-f}.
For simplicity we obtain the commutation relations between
currents $e^+(u)$ and $e^+(v)$ and set in the calculations $c=0$.
Others cases as well as the reconstraction of the central element
can be easily considered. We have
\bea
{[}e^+(u),e^+(v){]}&=& -i\h \int_{C_+}{d\tilde u\over 2\pi i}
\int_{C_+}{d\tilde v\over 2\pi i}\
{\{e(\tilde u),e(\tilde v)\}\over (\tilde u-\tilde v)(u-\tilde u)
(v-\tilde v)}=\nn\\
&=&
-i\h \int_{C_+}{d\tilde u\over 2\pi i}
\int_{C_+}{d\tilde v\over 2\pi i}
{e(\tilde u)e(\tilde v)\over \tilde u-\tilde v}
\left[{1\over (u-\tilde u)
(v-\tilde v)}-
{1\over (u-\tilde v)
(v-\tilde u)}\right]=\nn\\
&=&
i\h \int_{C_+}{d\tilde u\over 2\pi i}
\int_{C_+}{d\tilde v\over 2\pi i}
{e(\tilde u)e(\tilde v)\over u-v}
\left[{1\over u-\tilde u}-{1\over v-\tilde u}\right]
\left[{1\over u-\tilde v}-{1\over v-\tilde v}\right]=\nn\\
&=&i\h {(e^+(u)-e^+(v))^2\over u-v}\ ,\nn
\eea
where
the contour $C_+$ goes above the points
$u$, $v$ and in the second line  we used the fact that
$e^2(u)=0$ which follows from  the commutation relation \r{CC1}
in order to interchange the integration contours.
Due to this condition the pole of the integrand when
$\tilde u=\tilde v$ is superficial.
Let us stress that the condition that square of the total
current is vanishing is the special property of
deformed infinite dimensional algebras
(like central extended Yangian double, quantum affine algebra, etc.)
and has no classical analogs.

The next step is to obtain from the commutation relation
$$[e^+(u),e^+(v)]=i\h{(e^+(u)-e^+(v))^2\over u-v}\ $$
the relation \r{he-he} for $\la,\mu>0$. Others cases can be
obtained from the commutation relations $[e^\pm(u),e^\mp(v)]$
and $[e^-(u),e^-(v)]$. This can be done using
the convolution property of the
inverse Laplace transform and the following Lemma (recall that we consider
the case of $c=0$ for simplicity).
\begin{lemma}
$$
\hat g(\la,\mu)=
i\int_{C_+}{du\over 2\pi}\int_{C_+}{dv\over 2\pi}
\ee^{i\la u}\ee^{i\mu v}{e^+(u)-e^+(v)\over u-v}=
\left\{\begin{array}{ll}
\he_{\la+\mu}&\mbox{\rm for}\ \ \la+\mu>0,\\
\he_{\la}/2&\mbox{\rm for}\ \ \la>0,\ \mu=0,\\
\he_{\mu}/2&\mbox{\rm for}\ \ \mu>0,\ \la=0,\\
0&\mbox{\rm otherwise}.
\end{array}\right.
$$
\end{lemma}
Now we have
\bea
{[}\he_\la,\he_\mu{]}&=&
i\h
\int_{C_+}{du\over 2\pi}\int_{C_+}{dv\over 2\pi}
\ee^{i\la u}\ee^{i\mu v}\left[e^+(u) {e^+(u)-e^+(v)\over u-v}
- e^+(v) {e^+(u)-e^+(v)\over u-v}\right]
\nn\\
&=&\h \int_{-\infty}^\infty d\tau \int_{-\infty}^\infty d\nu
\ [\he_\tau \theta(\tau)\delta(\nu)-
\he_\nu \theta(\nu)\delta(\tau)]\ \hat g(\la-\tau,\mu-\nu)=\nn\\
&=&\h \int_{-\infty}^\infty d\tau\
[\theta(\tau-\mu)-\theta(\tau-\la)]\ \he_\tau \he_{\la+\mu-\tau}\ .
\label{result}
\eea
Note that the kernel of the integrand in \r{result} makes the
integral to be supported on the finite domain of the real axis.
One can verify that for $c\neq0$ the calculation \r{result} will
be modified a little bit but yields the same result \r{he-he}.

In the terms of the generators $\hhh_\la$ the commutation relations
\r{kk-hf} also can be written in a convolution form
:
\bea
{[}\hhh_\la,\he_\mu{]}&=&2\he_{\la+\mu} +\tih
\int_0^\la d\tau\ \left[\theta(\tau)-\theta(\tau-\la)\right]\
\{\hhh_{\tau},\he_{\la+\mu-\tau}\}\ ,
\label{hh-he}\\
{[}\tilde{\hhh}_\la,\hf_\mu{]}&=&-2\hf_{\la+\mu} -\tih
\int_0^\la d\tau\ \left[\theta(\tau)-\theta(\tau-\la)\right]\
\{\tilde{\hhh}_{\tau},\hf_{\la+\mu-\tau}\}\ ,
\label{hh-hf}
\eea
where $\tilde{\hhh}_\la=\hhh_\la\ee^{-\h c|\la|/2}$.

 The comultiplication rules \r{com-e-fu}--\r{comul-h} can be rewritten
for the generators $\he_\la$, $\hf_\la$ and $\hhh_\la$ in terms
of multiple convolution integrals. For example, at $c=0$ we have:
$$
\Delta\he_\la=\he_\la\ot 1 + 1\ot\he_\la +\h
\int^\la_0 d\tau\  \hhh_{\la-\tau}\ot \he_\tau
+o(\h^2).$$
The precise definition of the algebra $\Aya$ means that
one should consider the proper completion of the free
 tensor topological algebra generated by
$\he_\la$, $\hf_\la$, $\hhh_\la$,
 $\la \in\RR$, $c$ and
 $\tilde d$ over the ideal generated by the relations
 \r{hf-hf}.
 It means in particular that the completed algebra is generated by
 formal integrals
\beq
\int_{-\infty}^{+\infty}\he_\la g(\la)d\la ,\quad
 \int_{-\infty}^{+\infty}\hf_\la g'(\la)d\la ,\quad
 \int_{-\infty}^{+\infty}\kk_\la g''(\la)d\la ,
\label{425}
\eeq
where $g(\la)$, $g'(\la)$ and $g''(\la)$ are integrable functions decreasing
 faster then $\ee^{-a\mid\la\mid}$ for some positive $a$
 (see details in
 \cite{KLP3}) . The completion $\widehat{{\cal A}}$ is a Hopf algebra
  which coincides with rational degeneration $\Aya$ of the
 scaled elliptic algebra $\Ael$ \cite{KLP3}. As well as the Yangian double,
  it is the quantum double of its Hopf subalgebra, generated by
  the positive Fourier harmonics of the currents and one of the elements
  $c$ and $d$.
\bigskip

\noindent
{\bf \razdel The comparison.}
Let us compare the two Yangian type algebra
 obtained from the Yang $R$-matrix and the two Riemann problems. The
 commutation relations for their generators look very different while
 they were obtained from identical equations for the $L$-operators.
 In order to visualize their common nature one should look to the relations
 \r{comm2}, \r{comm4}--\r{comm7} as to the difference equations
 \r{comm4}--\r{comm7} with initial conditions \r{comm2}.
 One can solve these equations. For instance, the solution of \r{commut-dis}
for the generators $e_n$ has a form
 \beq
{[}e_k,e_l{]}=\yh\sum_{p=l}^{k-1}e_pe_{k+l-p-1}\ ,\qquad k>l
 \label{dis-conv}
\eeq
which is a discrete analog of the relation \r{he-he}.
Note that the relations \r{he-he},
\r{hf-hf}, \r{hh-he} and \r{hh-hf}
written in terms of
 convolutions allow to order the quadratic expressions over the generators
in terms of iterated integrals (see details in \cite{KLP3}).
On the other case, one can differentiate the relations \r{he-he},
\r{hf-hf}, \r{hh-he} and \r{hh-hf}
 over the parameters $\la$ and $\mu$ and get the continuous analogues of the
 relations \r{comm5}--\r{commut}:
\beq
{[}{\he'}_\la,{\he}_\mu{]}-{[}{\he}_\la,{\he'}_\mu{]}=
\h \{{\he}_\la,{\he}_\mu\}\ .
\label{diff1}
\eeq
 Nevertheless only the integral relations
\r{hf-hf}, \r{hh-he} and \r{hh-hf}
 have
 precise sense in completed algebra and admit further generalizations
 \cite{KLP3} to the contrary to their differential consequences \r{diff1}.

The natural question then arises why the algebras $\DY$ and $\Aya$ are
 essentially different. The main difference concerns the properties of the
 elements $d$ for $\DY$ and of $\tilde{d}$
  for the algebra $\Aya$ and the different types of completions.

In the discrete case of the Yangian double the operator $[d,\ \cdot \ ]
=\frac{d}{du}$ acts in $\DY$ as degree $-1$ operator preserving
 the filtration \r{2.7} and has no eigenvalues in the completed algebra
 different from zero. It would encounter the divergence of the partition
 functions of infinite-dimensional representations (see Section
\ref{section5}).

To the contrary the operator $[\tilde{d},\  \cdot \ ]
=\frac{1}{i}\frac{d}{du}$ act as degree zero operator in the algebra
 $\Aya$ preserving the grading. The elements $\he_\la, \hf_\la, \hhh_\la$
 are its eigenvectors with eigenvalue $\la, \la\in \RR$. As a consequence,
 the characters of integrable finite-dimensional modules are well-defined
 (see Section \ref{section5}).

The algebra $\Aya$ possesses the Cartan antiinvolution $\theta$ (which can
 be treated as an involutive antiisomorphism of the algebra $\Aya$ over the
 ring $\CC [\h]$ or as involutive antiisomorphism from $\Aya$ to
 ${\cal A}_{-\h}(\widehat{\frak sl}_2)$). On the level of generating
 functions it looks like
\beq
\theta e(u)=f(-u)\ , \qquad \theta f(u)=e(-u)
\ , \qquad \theta h^\pm(u)=\pm h^\mp(-u)\ ,
\label{theta}
\eeq
$$ \theta d =d\ , \qquad \theta c=c\  ,\qquad \theta \h=-\h\ ,$$
and for the generators
$$ \theta \he_\la =-\hf_{-\la},\qquad \theta \hf_\la =-\he_{-\la},\qquad
\theta \hhh_\la =\hhh_{-\la}\ .   $$
The antiinvolution $\theta$ exchange the two  subalgebras of $\Aya$
generated by Fourier modes of $e^\pm(u), f^\pm(u), h^\pm(u)$; thus the
 algebra $\Aya$ can be treated as contragredient algebra. To the contrary,
 there is no analog of Cartan antiinvolution for the Yangian double, since
 it reverses the filtration \r{2.7} and  should exchange Laurent  series
from a subalgebra generated by negative Fourier harmonics with Laurent
polynomials over positive Fourier harmonics. As a result we could not
 expect a good notion of restricted dual for highest weight modules
 of extended Yangian double.

The two algebras distinguish as well on the classical level.
The classical limit % ${\frak dy}(\widehat{{\frak sl}}_2)$
of the central  extended Yangian
double $\DY$
can be identify with semidirect sum of the central extension of
meromorphic $\frak{sl}_2$-valued functions and of one-dimensional Lie algebra
 generated by the derivative $\frac{d}{dz}$.
 The commutation relations for the generators are standard:
 $$ [h_n,e_m]=2e_{n+m},\qquad [h_n,f_m]=-2f_{n+m},\qquad
[e_n,e_m]=[f_n,f_m]=0,$$
 $$ [e_n,f_m]=h_{n+m}+n\delta_{n,-m}c,\qquad [h_n,h_m]=2n\delta_{n,-m}c$$
 $$ [d,e_n]=-ne_{n-1},\quad
 [d,f_n]=-nf_{n-1},\quad [d,h_n]=-nh_{n-1}\ .$$
 The bialgebra structure can be
defined via the decomposition of the algebra to the subalgebras of
 $\frak{sl}_2$-functions regular at zero and at infinity:
 $$ % {\frak dy}(\widehat{{\frak sl}}_2)=
\left({\frak sl}_2\ot{\CC}[z]\oplus{\CC}
 c \right)\oplus
 \left({\frak sl}_2\ot{\CC} [[z^{-1}]]\oplus{\CC} d\right)\ .$$
 Again we have no symmetry between the two subalgebras and
 the (quasi)nilpotent action   of the operator $d$.

 The classical limit of the Hopf algebra $\Aya$ can be
  identify with the completion of central extended algebra of
 meromorphic  $\frak{sl}_2$-valued functions vanishing at infinity
 (again with added element $\tilde{d}= \frac{1}{i}\frac{d}{dz}$).
 The commutation relations for the generators now are
 $$ [\hhh_\la,\he_\mu]=2\he_{\la+\mu},\qquad [\hhh_\la,\hf_\mu]
=-2\hf_{\la+\mu},\qquad [\he_\la,\he_\mu]=[\hf_\la,\hf_\mu]=0,$$
 $$ [\he_\la,\hf_\mu]=\hhh_{\la+\mu}+\la\delta(\la+\mu)c,\qquad
  [\hhh_\la,\hhh_\mu]=2\la\delta(\la+\mu)c,$$
 $$
 [\tilde{d},\he_\la]=-\la \he_{\la},\quad
 [\tilde{d},\hhh_\la]=-\la \hhh_{\la},\quad
 [\tilde{d},\hf_\la]=-\la \hf_{\la}.$$
  The bialgebra structure is given by the decomposition of
 the algebra into direct sum of subalgebras of
 $\frak{sl}_2$-valued functions which are regular at upper or
  lower half-plane, which corresponds to taking positive or
  negative Fourier modes. Again we have the contragredient structure
  and the diagonalization of the operator $\tilde{d}$.

 Note also that the different types of completions used for
 the definitions of the two algebras actually follow from different
 positions of the marked singular points in two Riemann problems.
 In the case of the
 Yangian double two separated points zero and infinity define
 Laurent series in infinity and Laurent polynomial in zero while
in the case of $\Aya$ the singular point near the contour (infinity)
produces continuous family of Fourier harmonics which diagonalize the
 operator $\tilde{d}$. Thus these two algebras cannot be connected
 by a projective transform of the complex plane as well as the
 corresponding Riemann problems.

Let us note that the algebras $\DY$ and $\Aya$ given by the commutation
relations \r{comm5}--\r{commut} and \r{kk-hf}--\r{he-hf}
respectively and associated with simple Lie algebra $\frak{sl}_2$
can be generalized for arbitrary simply-laced Lie algbera.
See, for example, \cite{K}.

\setcounter{equation}{0}
\setcounter{parag}{0}

\section{Finite-dimensional representations and $R$-matrices}
{\bf \razdel Finite-dimensional representations.}
The finite-dimensional representations make sense for the
 subquotions of the algebras $\DY$ and $\Aya$ defined by the
 condition $c=0$ with dropped elements $d$ and $\tilde{d}$.
 Thus their structures are identical for both
 algebras. Due to classification theorem (see \cite{Chary}),
the irreducible finite-dimensional representations of the
 algebras $\DY$ and $\Aya$ are certain tensor products of
 the evaluation representations. The evaluation representations
 can be defined via evaluation homomorphism $\Ev_z :\DY\to U(\frak{sl}_2)$,
 $\Aya\to U(\frak{sl}_2)$, $z\in\CC$. For the Yangian double the homomorphism
  has a form
\beq
\Ev_z \sk{e(u)}=\yh
\delta( u-z-\frac{h-1}{2}\yh)\cdot e\ , \quad
\Ev_z \sk{f(u)}=\yh
 f\cdot\delta( u-z-\frac{h-1}{2}\yh)\ ,
\label{ev}
\eeq
\beq
\Ev_z \sk{h^\pm(u)}=1 + {\yh \over
u-z-\yh(h-1)/2}\, ef -{\yh \over u-z-\yh(h+1)/2} fe\ ,
\label{eval-h}
\eeq
as well as for the algebra $\Aya$ with the replacement $\yh=-i\h$.
 Here $e, f$   and $h$ are the generators of Lie algebra  $\frak{sl}_2$:
$$ [h,e]=2e,\qquad [h,f]=-2f,\qquad [e,f]=h $$
 and the right hand sides of \r{eval-h} are taken for
 $\Im z \gtrless \Im u$ for $\Aya$ and $|z|\gtrless|u|$ for $\DY$
respectively for ``$\pm$'' generating functions.
For instance, for the simplest two-dimensional representation
 $\pi_z^{(1/2)}$ the action of the generators
 of the algebra $\DY$ and of the algebra $\Aya$ has a form
\beq
e_k= z^k e_{1,2},\quad f_k=z^k e_{2,1},\quad h_k= z^k (e_{1,1}-e_{2,2})\ ,
\label{twodim1}
\eeq
\bea
\he_\la&=&\ee^{i\la z}e_{1,2},\quad
\hf_\la=\ee^{i\la z}e_{2,1},\quad
\hhh_\la=\ee^{i\la z}(e_{1,1}-e_{2,2})\ ,\nn\\
\kk_\la&=&{\ee^{i\la z}\over \h\la}
\left((1-\ee^{-\la\h})e_{1,1}+(1-\ee^{\la\h})e_{2,2}\right),
\label{twodim2}
\eea
 where $(e_{i,j})_{k,l}=\delta_{ik}\delta_{jl}$
are unit matrices in $\mbox{End}\ \CC^2\ot\CC^2$.
\bigskip

\noindent
{\bf \razdel Universal ${\cal R}$-matrices.}
The universal ${\cal R}$-matrix for the algebra $\DY$
has been obtained in \cite{KT,K} from the analysis of
the canonical Hopf pairing  of the two Hopf subalgebras $\DY^\pm$
 of $\DY$, generated by Fourier coefficients of the currents
 $e^\pm(u), f^\pm(u)$, $h^\pm(u)$.
The pairing of these fields looks like
\beq
\langle e^+(u),f^-(v)\rangle =\langle f^+(u),e^-(v)\rangle =
\frac{\yh}{u-v}\ ,
\label{pairing1}
\eeq
\beq
\langle h^+(u),h^-(v)\rangle =\frac{u-v+\yh}{u-v-\yh}
\label{pairing2}
\eeq
or, in terms of the generators,
$$%\beq
\langle c,d\rangle=-\yh^{-1},\quad
\langle e_k,f_{-l-1}\rangle=\langle f_k,e_{-l-1}\rangle
=-\yh^{-1}\delta_{k,l}\ ,
$$
$$
 \langle h_k, h_{-l-1}\rangle =
\left\{\begin{array}{l}
2(\yh)^{k-l-1}{k!\over l!(k-l)!},\quad
k\geq0,\quad  0\leq l\leq k\ ,\\
0,\quad \hbox{otherwise}\ .\end{array}\right.
$$%\label{pair-dis}\eeq
The full pairing is described by the universal
$R$-matrix which has a form
\beq
{\cal R} = {\cal R}_+\cdot {\cal C}\cdot
{\cal R}_0 \cdot {\cal C} \cdot {\cal R}_-
\ ,
\label{univ-R}
\eeq
where
\bea
{\cal C}&=&-{\yh\over 4} (c\ot d+ d\ot c)\ ,\nn\\
{\cal R}_+&=& \prod_{k\geq 0}^\rightarrow \exp (-\yh e_k\ot f_{-k-1})
= \exp(-\yh e_0\ot f_{-1})\exp(-\yh e_1\ot f_{-2})\cdots
\ ,\label{R+}
\\
{\cal R}_-&=& \prod_{k\geq 0}^\leftarrow \exp (-\yh f_k\ot e_{-k-1})
= \cdots \exp(-\yh f_1\ot e_{-2})\exp(-\yh f_0\ot e_{-1})
\ ,\label{R-}\\
{\cal R}_0&=& \prod_{n\geq 0} \exp
\left(-\Res{u=v}
\left[ {d\over du} \ln\, h^+(u)\ot
\ln\,h^-(v-(2n+1)\yh)\right]\right)
\ ,\label{R-01}
\eea
and a residue operation $\Res{}$ is defined as follows
$$
\Res{u=v}\left(\sum_{i\geq0} a_i u^{-i-1}\ot
\sum_{k\geq0} b_k v^{k}
\right) =
\sum_{i\geq0} a_i\ot b_i\ .
$$
As usual, the $L$-operators $L^\pm$ are given by the substitution to
one tensor component of ${\cal R}$ the two-dimensional representation
\cite{FR}
\beq
L^-(z)=(\pi^{(1/2)}_z\ot\id)
\ {\cal C}^{-1}\cdot{\cal R}\cdot{\cal C}^{-1},\quad
L^+(z)=(\pi_z^{(1/2)}\ot\id)\
{\cal C}\cdot\left({\cal R}^{21}\right)^{-1}\cdot{\cal C}\ .
\label{L-operat}
\eeq
The decomposition \r{univ-R} produces the Gauss decomposition of the
$L$-operators \r{GL-univ}.

Analogously, the description of the universal $R$-matrix for the algebra
$\Aya$ can be done following the same scheme of \cite{KT}. Let us remind
that the main arguments in \cite{KT} are: the triangular decomposition
of the Yangian Double with respect to a Hopf pairing, the basic pairing
\r{pairing1}, \r{pairing2}, and the expression of the tensor of the pairing for
the subalgebras generated by Cartan currents as the exponent of the
 pairing between logarithms of Cartan fields. The last calculation uses
 the shift automorphisms of $\DY$.
 All these arguments remain
 unchanged for the algebra $\Aya$ (in basic pairing \r{pairing1},
 \r{pairing2} we use $\yh = -i\h$ as usual).

  It means that the universal $R$-matrix for $\Aya$ admits the decomposition
 \r{univ-R} where the factor ${\cal C}$ is unchanged, the ordered
 products of exponents in the factors ${\cal R}_\pm$ turn to ordered
 exponential in integral form and the factor ${\cal R}_0$ is as before
 $${\cal R}_0 =\exp \Omega$$
 where $\Omega$ is the tensor of the pairing of the fields $k^\pm(u)=
\log h^\pm(u)$. The main distinction of the case of $\Aya$ is that the
 pairing
 $$\langle \kappa^+(u),\kappa^-(v)\rangle=\log \frac{u-v-i\h}{u-v+i\h}$$
 can be explicitly and uniquely diagonalized in the generators $\kk_\la$:
 $$\langle \kk_\la,\kk_\mu\rangle
=-2\frac{\mbox{sh}\h\la}{\h^2\la}\delta(\la+\mu)\ .   $$
Summarizing the calculation we have.

The universal $R$-matrix for $\Aya$ admits the decomposition \r{univ-R}
 where
${\cal C}=\exp(-\h(\tilde d\ot c+c\ot \tilde d)/4)$,
\beq
%{\cal R}_0&=& \exp \left(\int_0^\infty {d\la \h\la\over 2\sh\,\h\la}
%\hhh_{\la}\ot \hhh_{-\la}\right)\ ,\nn\\
{\cal R}_+\ =\  \mathop{P}\limits^\rightarrow \exp
\left(-\h\int_0^{+\infty} d\la\,
\he_{\la}\ot \hf_{-\la}\right),\quad %\nn\\
{\cal R}_-\ =\  \mathop{P}\limits^\leftarrow \exp
\left(-\h\int_0^{+\infty} d\la\,
\hf_{\la}\ot \he_{-\la}\right)\ .
\label{Rpm}
\eeq
and
\beq
{\cal R}_0= \exp \left(-\int_0^{+\infty} d\la\  {\h^2\la\over 2\,\sh\,\h\la}
\kk_{\la}\ot \kk_{-\la}\right)\ .
\label{R00}
\eeq

It is interesting to compare the evaluation of the two expressions
 of the universal $R$-matrix to tensor product $\pi_{z_1}^{(1/2)}\ot
\pi_{z_2}^{(1/2)}$ of two-dimensional representations.
The formulas \r{R+}, \r{R-} and \r{R-01} yield the following
 decomposition of four by four $R$-matrix $R(z)$:
\beq
R^-(z)=\left(
\begin{array}{cccc}
1&0&0&0\\
0&1&\frac{\yh}{z}&0\\
0&0&1&0\\
0&0&0&1
\end{array}
\right)
\left(
\begin{array}{cccc}
\rho^-(z)&0&0&0\\
0&\frac{z-\yh}{z}\rho^-(z)&0&0\\
0&0&\frac{z}{z+\yh}\rho^-(z)&0\\
0&0&0&\rho^-(z)
\end{array}
\right)
\left(
\begin{array}{cccc}
1&0&0&0\\
0&1&0&0\\
0&\frac{\yh}{z}&1&0\\
0&0&0&1
\end{array}
\right)
\label{appl}
\eeq
or $R^-(z)=\rho^-(z)\frac{z+\yh P}{z+\yh}$,
where $P$ is a permutation of tensor
components and
$$\rho^-(z)=\prod_{n\geq 0}\frac{(z-2n\yh)(z-(2n+2)\yh)}{(z-(2n+1)\yh)^2}\ .$$
The infinite product converge for proper $z$ and equals to the ratio
 of $\Gamma$--functions:
$\rho^-(z)=\frac{\Gamma^2({1\over 2}-{z\over 2\yh})}
{\Gamma(1-{z\over 2\yh})
\Gamma(-{z\over 2\yh})}$.
An application of \r{Rpm} and \r{R00} gives $R^-(z)$ in a form \r{appl}
 but with a scalar factor $\rho^-(z)$ presented in an integral form
$$\rho^-(z)=\exp \left(
-2\int_0^{+\infty} d\la\ \frac{\sh^2({\h\la}/{2})}{\la\,\sh\, \h\la}
\ee^{i\la z}\right),\quad \Im\, z>0\ .$$
We see that the two universal $R$-matrices presented in this section
give two different quantization of the Yang rational solution of classical
 YB equation. The solution via algebra $\Aya$ has an advantage of
 integral presentation and uniquely determined by the asymptotics of
 the scalar factor. To the contrary, there is no definite choice
 for such a solution for the double of the Yangian (the answer has no
 definite asymptotics when $|z|$ tends to infinity). Moreover,
 as we will see further analogous integral representation
 for divergent infinite products automatically appear
 in representation theory of $\Aya$ in regularized form whereas
 the representation theory of the Yangian double has no instruments
 for such a regularization.

\setcounter{equation}{0}
\setcounter{parag}{0}

\label{section5}
\section{Basic infinite-dimensional representations}

The goal of this section is to construct the examples of the
infinite-dimensional rep\-re\-sen\-ta\-tions of the algebras $\DY$ and
$\Aya$. This construction demonstrate the distinctions of the
considered algebras. We start from the discrete algebra
$\DY$ and for simplicity consider only the case of the central
element $c=1$.
\medskip

\noindent
{\bf \razdel
Basic representations of the algebra $\DY$.}
Let ${\cal H}$ be the Heisenberg algebra generated by free
bosons $a_{\pm n}$, $n=1,2,\ldots$ with zero modes
$a_0$, $\da$  and commutation relations
$$%\begin{equation}
%\label{bosons}
[a_n,a_m] =  n \delta_{n+m,0}\ ,\quad [\da,a_0]=2\ .
$$%\end{equation}
Let
$$
a_+(z)=\sum_{n\geq 1}\frac{a_n}{n}z^{-n}-p\log z\ ,\quad
a_-(z)=\sum_{n\geq 1}\frac{a_{-n}}{n}z^{n}+\frac{a_0}{2}\ ,%\label{5.2}
\quad
\phi_\pm(z)=\exp a_\pm(z)\ ,
$$
be the generating functions of the elements of the algebra ${\cal H}$.
Let $\bar e(u)$, $\bar f(u)$, $\bar h^+(u)$ and $\bar h^-(u)$
be following generating series acting in the Fock space $\HH$:
\bea
\bar e(\u)&=&\yh\phi_-(\u-\yh)\phi_-(\u)\phi^{-1}_+(\u)\ ,\nn\\
\bar f(\u)&=&\yh\phi^{-1}_-(\u+\yh)\phi^{-1}_-(\u)\phi_+(\u)\ ,\nn\\
\bar h^+(\u)&=&\phi_+(\u-\yh)\phi_+^{-1}(\u)\ ,\quad
\bar h^-(\u)=\phi_-(\u-\yh)
\phi_-^{-1}(\u+\yh)\ .\nn
\eea
We have the following

\begin{prop}\label{bos-dis}
{\rm\cite{K}}
{\sl The ${\cal H}$-valued generating functions (fields)}
$$
e(u)=\bar e(u+c\yh/4),\quad
f(u)=\bar f(u-c\yh/4),\quad
h^+(u)=\bar h^+(u+c\yh/2),\quad
h^-(u)=\bar h^-(u)
$$
{\sl satisfy the commutation relations {\rm\r{CCeta}}
with $c=1$.}
\end{prop}

Let $V_\alpha$ be the formal power series extensions of the Fock spaces
\beq
V_i={\CC}[[a_{-1},\ldots ,a_{-n},\ldots ]]\otimes
        \left(\oplus_{n\in{\ZZ}}{\CC}\ee^{(n+\alpha)a_0}\right)\ ,
\quad 0\leq \alpha < 1\ ,
\label{modules}\eeq
with the action of bosons on these spaces
\bea
a_n&=&\hbox{the left multiplication by $a_n\otimes1$\ \  for $n<0$}\ ,\nn\\
   &=&{[}a_n,\ \cdot\ {]}\otimes 1\quad \hbox{for $n>0$}\ ,\nn\\
\ee^{n_1a_0}(a_{-j_k}\cdots a_{-j_1}\otimes\ee^{n_2a_0})&=&
a_{-j_k}\cdots a_{-j_1}\otimes\ee^{(n_1+n_2)a_0}\ ,\nn\\
u^{\da} (a_{-j_k}\cdots a_{-j_1}\otimes\ee^{na_0})&=&u^{2n}
a_{-j_k}\cdots a_{-j_1}\otimes\ee^{na_0}\ .\label{action}
\eea
It is clear from \r{action} that Fock spaces $V_\alpha$ becomes
the irreducible representations of the algebra $\DY$ at level 1
($c=1$) for $\alpha=0$ or $1/2$ with the vacuum vectors $1\ot1$
and $1\ot \ee^{a_0/2}$. These representations are highest
weight representations with respect to the Fourier components
$\bar e_n$, $\bar f_n$
of the generating currents $\bar e(u)=\sum \bar e_n u^{-n-1}$
and $\bar f(u)=\sum \bar f_n u^{-n-1}$. The latter generators  are related to
$e_n$ and $f_n$ by the triangular transformations
due to the relations given in Proposition~\ref{bos-dis}
\bea
\bar e_n(1\ot1)&=&\bar f_n(1\ot1)=0, \quad n<0,\nn\\
\bar e_n(1\ot\ee^{a_0/2})&=&0,
\quad n<-1, \quad
\bar f_n(1\ot\ee^{a_0/2})=0,
\quad n<1.\nn
\eea

Elements of the monomial basis
\r{modules} are {\em not} eigenvalues of the filtration operator $d$,
since:
\bea
{[}d,a_n{]}&=&-na_{n-1},\quad n\leq -1, \quad n\geq2\ ,\nn\\
{[}d,a_1{]}&=&-p,\quad {[}d,p{]}=0,\quad {[}d,a_0/2{]}=a_{-1}\ .\nn
\eea
They cannot be used for example for calculation the character of the
Fock space $V_\alpha$ using the operator $d$. Usual expression
$\tr_{V_\alpha}(\ee^{pd})$ is divergent and only the ratio of
such traces
$\left.\tr_{V_\alpha}(\ee^{pd}O)\right/\tr_{V_\alpha}(\ee^{pd})$
can be made finite for certain operators $O$ in the Fock
space $\HH$  \cite{C}.

\noindent
{\bf \razdel Realization  of the algebra $\Aya$
by the continuous fields.}
Let us define bosons $a_\la$, $\la\in\RR$ which satisfy
the commutation relations \cite{L,JKM}:
\beq
[a_\la,a_\mu]={4\over\tih^2}\,{\sh(\tih\la)\sh(\tih\la/2)\over\la}\,
\, \delta(\la+\mu)=
a(\la)\delta(\la+\mu)\ .
\label{bosons}
\eeq
Consider the generating functions
\bea
e(u)&=&\h\ee^\gamma{:}\exp\left(-\tih \intt d\la\ \ee^{-i\la u}{a_\la
\ee^{-\h|\la|/4}\over
2\sh(\tih\la/2)}\right){:}\ ,\label{Ebos}\\
f(u)&=&\h\ee^\gamma{:}\exp\left(\tih \intt d\la\ \ee^{-i\la u}{a_\la
\ee^{\h|\la|/4}\over
2\sh(\tih\la/2)}\right){:}\ ,\label{Fbos}\\
h^\pm(u)&=&(\h\ee^{\gamma})^{-2}
\ {:}\ e\left(u\mp\fract{i\tih}{4}\right)
f\left(u\pm\fract{i\tih}{4}\right){:} =\nn\\
&=& \exp\left(\tih \int_0^{\pm\infty} d\la\ \ee^{-i\la u}
a_\la\right)\ ,
\label{Hbos}
\eea
where $\gamma$ is Euler constant.
The notion of the normal ordered operator becomes more involved in case of
continuous bosons. It requires some kind of ultraviolet regularization
to be included in the definition of the normal ordered operators such
that the product of these operators satisfy the rules
\cite{JM}:
\bea
&{:}\exp\left(\intt d\la\ g_1(\la)\,a_\la\right){:}\ \cdot\
{:}\exp\left(\intt d\mu\ g_2(\mu)\,a_\mu\right){:} = \nn\\
&\kern-30pt =\exp\left(\int_{\tilde C}{d\la\,\ln(-\la)\over2\pi i}\
a(\la)g_1(\la)g_2(-\la) \right)
{:}\exp\left(\intt d\la\ (g_1(\la)+g_2(\la))\,a_\la\right){:}\ .
\label{normal}
\eea
The contour $\tilde C$ is shown in the Fig.~3.
\bigskip

\unitlength 1.00mm
\linethickness{0.4pt}
\begin{picture}(121.00,20.00)
\put(17.00,15.00){\makebox(0,0)[cc]{0}}
\put(20.00,15.00){\makebox(0,0)[cc]{$\bullet$}}
\put(132.00,15.00){\makebox(0,0)[cc]{$+\infty$}}
\put(20.00,15.00){\line(1,0){100.33}}
\put(40.00,10.00){\line(1,0){80.33}}
\put(120.00,20.00){\line(-1,0){100.00}}
\put(30.00,5.00){\makebox(0,0)[cc]{Fig.~3.}}
\put(121.00,10.00){\vector(1,0){0.2}}
\put(100.00,10.00){\line(1,0){21.00}}
\put(20.00,10.00){\line(1,0){22.00}}
\put(20.67,10.00){\line(1,0){22.00}}
%\put(42.67,10.00){\line(0,1){0.00}}
\put(20.17,15.00){\oval(15.00,10.00)[l]}
\end{picture}
\smallskip

With this definition of the normal ordered exponents we can prove the
following
\begin{prop}\label{repr-con}
{\sl The generating functions} {\rm\r{Ebos}--\r{Hbos}}
{\sl satisfy the commutation relations  {\rm\r{CC1}}}.
\end{prop}

The main formula which should be used to prove the Proposition
\ref{repr-con} is:
$$
\exp\left(
\int_{\tilde C}{d\la\,\ln(-\la)\over 2\pi i\la}\ee^{-x\la}
\right)=\ee^{-\gamma}
x^{-1}, \quad \Re x>0\ .
$$

For the description of the infinite-dimensional
representations of the algebra $\Aya$ in terms of continuous bosons
we need a definition of a Fock space generated by the continuous
family of free bosons. We borrow the construction below from
\cite{KLP3}.

Let $a(\la)$ be a meromorphic function, regular for $\la\in \RR$ and satisfying
 the following conditions:
$$a(\la)= -a(-\la)\ ,$$
$$
a(\la) \sim a_0 \la,\quad \la \to 0, \qquad a(\la)\sim
\ee^{a' |\la|},  \quad \la \to \pm\infty .
$$

Let $a_\la$, $\la\in\RR$, $\la\neq 0$ be free bosons which satisfy the
 commutation relations
$$%\beq
[a_\la, a_\mu]=a(\la)\delta(\la+\mu) \ .
$$%\label{boson}\eeq

 We define a  (right) Fock space $\HH_{a(\la )}$ as follows.
$\HH_{a(\la )}$ is generated as a vector
space
by the expressions
$$
\mint f_n(\la_n) a_{\la_n} d\la_n\ldots \mint f_1(\la_1) a_{\la_1} d\la_1\
\rvac\  , $$
 where the functions $f_i(\la)$ satisfy the condition
$$%\begin{equation}
f_i(\la)< C\, \ee^{(a'/2+\epsilon)\la}, \qquad \la \to -\infty\ ,
$$%\end{equation}
for some $\epsilon>0$ and
$f_i(\la)$ are analytical functions in a neighborhood of
 ${\RR_-}$ except $\la=0$, where they have a simple pole.

The left Fock space $\HH^*_{a(\la)}$ is generated by the expressions
$$
\lvac \nint g_1(\la_1) a_{\la_1} d\la_1\ldots \nint g_n(\la_n) a_{\la_n}
d\la_n\ ,
$$
where the functions $g_i(\la)$ satisfy the conditions
$$%\begin{equation}
g_i(\la)< C\, \ee^{-(a'/2+\epsilon)\la}, \qquad \la \to +\infty\ ,
$$%\end{equation}
for some $\epsilon>0$
and $g_i(\la)$ are analytical functions in a neighborhood
of
 ${\RR_+}$ except $\la=0$, where they also have a simple pole.

The pairing $(,):$ $\HH_{a(\la )}^*\otimes \HH_{a(\la )} \to \CC$ is
uniquely defined by
the
 following prescriptions:
\bea
 &(i)&\ \  (\langle\mbox{vac}|, \rvac) =1\ ,\nn\\
 &(ii)&\ \  (\lvac \nint d\la \ g(\la)a_\la \ ,
\mint d\mu\  f(\mu)a_\mu\   \rvac) =
\int_{\tilde C} {d\la\,\ln(-\la)\over 2\pi i}
  g(\la)f(-\la)a(\la)\ ,\nn\\
&(iii)&\ \ \mbox{the Wick theorem}.\nn
\eea

Let the vacuums $\lvac$ and $\rvac$ satisfy
 the conditions
$$
a_\la\rvac =0,\quad \la>0,\qquad \lvac a_\la =0, \quad \la<0\ ,
$$
 and $f(\la)$ be  a function analytical in some neighborhood of the real line
 with possible simple pole at $\la=0$ and which has the following asymptotical
 behavior:
$$
\qquad f(\la)< C e^{-(a'/2+\epsilon)|\la|}
, \qquad \la \to \pm\infty
$$
for some $\epsilon>0$.
Then, by definition, the operator
$$F= {:}\exp \sk{\int_{-\infty}^{+\infty}d\la\ f(\la)a_\la} {:}$$
acts on the right Fock space $\HH_{a(\la )}$ as follows.
$F=F_-F_+$, where
$$F_-=\exp \sk{\int_{-\infty}^{0}d\la\ f(\la)a_\la} \quad
\mbox{and}\quad
F_+=
\lim{\epsilon\to +0} \ee^{
\epsilon \ln \epsilon f(\epsilon)a_{\epsilon}}
\exp\sk{\int_{\epsilon}^{\infty}d\la\ f(\la)a_\la}\ .
$$
The action of operator $F$ on the left Fock space $\HH_{a(\la )}^*$ is
 defined via another decomposition: $F=\tilde{F}_-\tilde{F}_+$, where
$$\tilde{F}_+=\exp \sk{\nint d\la\  f(\la)a_\la}\quad
\mbox{and}\quad
\tilde{F}_-=
\lim{\epsilon\to +0} \ee^{
\epsilon \ln\epsilon f(-\epsilon)a_{-\epsilon}}
\exp\sk{\int^{-\epsilon}_{-\infty}d\la\ f(\la)a_\la}.
$$
 These definitions imply the following statement:
\smallskip

\begin{prop}\label{777}

\begin{itemize}
\item[$(i)$]{\sl The defined above actions of the operator
$$
F= :\exp \sk{\int_{-\infty}^{+\infty}d\la\ f(\la)a_\la}{:}
$$
on the Fock spaces $\HH$ and $\HH^*$ are adjoint};
\item[$(ii)$]{\sl The
product of the normally ordered operators satisfy the property}
 {\rm\r{normal}}.
\end{itemize}
\end{prop}

Returning to level one representation of $\Aya$ we choose $\HH =
\HH_{a(\la )}$
 for $a(\la)$ defined in \r{bosons}:
$$a(\la)={4\over\tih^2}\,{\sh(\tih\la)\sh(\tih\la/2)\over\la}\ .$$
From the definition of the Fock space $\HH$ and from the proposition 6 we have
 immediately the construction of a representation of $\Aya$:
\smallskip

\begin{prop}
\label{8}
{\sl The relations {\rm\r{Ebos}--\r{Hbos}} define a highest
 weight right representation of the algebra $\Aya$
in the Fock space $\HH$ and lowest weight left representation in the
dual Fock space $\HH^*$:}
\end{prop}
\beq
\he_\la\vac=0,\quad \hf_\la\vac=0,\quad \la\geq0\quad\mbox{and}\quad
\lvac\he_\la=0,\quad \lvac\hf_\la=0,\quad \la\leq0\ .
\label{hwr}
\eeq

The highest weight property \r{hwr} means that all the matrix elements
of the corresponding operators which do not vanishing identically
satisfy this property. Let us demonstrate that
$\langle v\mid \he_\la\vac =0$ for  $\la>0$ and certain $v\in\HH^*$.
Fix $\h>0$.
It is clear that any such matrix element has a form:
\bea
&\lvac \prod_{i=1}^n f(v_i)\prod_{j=1}^{n-1} e(u_j) \he_\la\vac=\nn\\
&\qquad ={1\over 2\pi}\int_{-\infty}^{\infty}du_n\ \ee^{i\la u_n}
\lvac \prod_{i=1}^n f(v_i)\prod_{j=1}^{n} e(u_j) \vac=\nn\\
&\qquad = G(v;u)\int_{-\infty}^\infty
du_n\ \ee^{i\la u_n}{\dis \prod_{j=1}^{n-1}(u_n-u_j)(u_n-u_j+i\h)
\over \dis \prod_{i=1}^{n}(u_n-v_i+i\h/2)(u_n-v_i-i\h/2)}\ ,
\label{cal11}
\eea
where $G(v;u)$ is some factor which does not depend on the variable
$u_n$ and in order to calculate \r{cal11} using
normal ordering relations for
total currents (see \cite{KLP3} for details)
 we should impose $\Im\, v_i<-\h/2$,
$i=1,\dots,n$ and $\Im\, u_j<0$, $j=1,\ldots,n-1$. Note that
integrand is decreasing as $u_n^{-2}$ function when $u_n\to\pm\infty$
so the integral is convergent. To calculate it we can close the contour
of integration either along the big semicircle in upper half-plane of
$u_n$ for $\la>0$ or along big semicircle in lower half-plane.
But all the poles of the integrand are in lower half-plane so
for $\la>0$ the nonvanishing matrix elements of the operator $\he_\la$
equal to zero and the property $\he_\la\vac=0$ for $\la>0$ is proved.
When $\la=0$ \r{hwr} follows from the continuity  arguments.

In the contrary to the case of discrete algebra $\DY$ the trace
function $\tr_{\cal H}(\ee^{pd})$ is well defined now and
can be calculated using the gradation property of the operator
$d$:
$$
[d,a_\la]=-\la a_\la\ .
$$
By definition this trace is equal to:
\bea
\tr_{\cal H}(\ee^{pd})&=&\sum_{n=0}^\infty\mathop{\int\cdots\int}
\limits_{0\leq\la_1<\cdots<\la_n<\infty}
d\la_1\cdots d\la_n
{\lvac a_{\la_1}\cdots a_{\la_n} \ee^{pd} a_{-\la_n}\cdots a_{-\la_1} \rvac
\over
\lvac a_{\la_1}\cdots a_{\la_n} a_{-\la_n}\cdots a_{-\la_1} \rvac}
\nn\\
&=&\sum_{n=0}^\infty {1\over n!} \left(\int_{0}^\infty d\la\ \ee^{p\la}
\right)^n= \ee^{1/p},\quad \Re p <0\ .\label{trace-res}
\eea
For the interpretation the generalized form-factors in quantum
integrable field theory as traces \r{trace-res}
one should put usually
$\Re p=0$, ($p=2\pi i$). In this case we should understand the result
\r{trace-res} as analytically continued from the domain, where $\Re p<0$.
This result can be compared with asymptotical expansion of the
partition function $\prod_{n=1}^{\infty}(1-q^n)^{-1}$. Really,  the trace
of the operator $\ee^{pd}$ can be presented as the integral over
eigenvalues of this operator $\ee^{p\la}$:
$$
\ee^{1/p}=1+\int_{0}^\infty d\la\ \ee^{p\la} \frak{p}(\la)\ ,
$$
where
\beq
\frak{p}(\la)=\sum_{n\geq0}{1\over n!(n+1)!} \la^n\ .
\label{partition-con}
\eeq
On the other hand the coefficients $\frak{p}_n$ of the partition function
$$
\prod_{n>0}{1\over (1-q^n)}=\sum_{n\geq 0} \frak{p}_n q^n
$$
have following asymptotical expansion
\beq
\frak{p}_n\sim\sum_k {1\over k!(k-1)!} n^k
\label{partition-dis}
\eeq
in the region when $n\to\infty$.
The explicit comparing \r{partition-con} and \r{partition-dis}
demonstrated their similarity.

\setcounter{equation}{0}
\setcounter{parag}{0}
\section{Quantized Current Algebras}

{\bf \razdel}
As we have seen in the previous section the total current algebra is
apparently more suitable for the constructing the infinite-dimensional
representations than standard $RLL$-formalism
which used the Gauss coordinates
of $L$-operators. On the other hand
it is difficult to define the Hopf structure in terms of total
currents  while in the $L$-operator formalism the Hopf structure
$\Delta' L= L\ot L$ is quite natural due to $RLL$-relations and
YB equation for $R$-matrix. We would like to discuss in this
section the following point of view.  It is possible to assign the
Hopf algebra structure to the algebra of total currents adding to the
commutation relations \r{CC1} the information on the Riemann problem
and also some additional information which we will discuss below.

Let us start with formal total current algebra for the
currents $e(u)$, $f(u)$ and $h^\pm(u)$ given for example
by the commutation relations \r{CC1}. This algebra is formal
since in the commutation relations of total currents $[e(u),f(v)]$
the $\delta$-function is defined formally as well as  Cartan "half"
currents $h^\pm(u)$.

The first step is to fix the Riemann problem which serves to divide
the total currents $e(u)$ and $f(u)$ into "half"-currents $e^\pm(u)$
and $f^\pm(u)$. This fixes the notion of $\delta$-function in the
commutation relations $[e(u),f(v)]$, the contents of the
generating functions $h^\pm(u)$
and also the full set of the commutation relations between all
the generating functions $e^\pm(u)$, $f^\pm(u)$ and $h^\pm(u)$
(see the Proposition \ref{3}). On the other hand fixing first the
analytical properties of $\delta$-functions in this relation leaves
the freedom for the factorization problem. This freedom is related
to the possible twist in the Riemann problem and will be also
discussed below.

Let us assume now that
the generating functions
$e^\pm(u)$, $f^\pm(u)$ and $h^\pm(u)$ are Gauss coordinates
of some $L$-operator given by \r{GL-univ}.
It is natural
to guess that this $L$-operators satisfy $RLL$-relations with some
$R$-matrix.

We state that there exist the universal comultiplication rules
for Gauss coordinates
$e^\pm(u)$, $f^\pm(u)$,
$h^\pm(u)=k^\pm_1(u)\left(k^\pm_2(u)\right)^{-1}$
and
$\tilde h^\pm(u)=\left(k^\pm_2(u)\right)^{-1}k^\pm_1(u)$
($h^\pm(u)$ may be not equal to
            $\tilde h^\pm(u)$  in the   general case)
which are based on the only assumption that the $R$-matrix
in a critical point is given by the rank one operator.
Let, for instance, $\overline R(i\h)$ is proportional to
$1-P$, where $P$ is a flip:
$$
\overline R(i\h)=\left(\begin{array}{cccc}
0&0&0&0\\
0&1&-1&0\\
0&-1&1&0\\
0&0&0&0
\end{array}\right).
$$
It takes place for the Yang $R$-matrix, for Baxter's elliptic $R$-matrix and
for Sine-Gordon $R$-matrix.  Then
the commutation relations for the $L$ operators at the
critical point of $R$-matrix basically reduces to the following:
\bea
k^\pm_1(u) f^\pm(u+i\h)&=&f^\pm(u) k^\pm_1(u)\ ,\nn\\
k^\pm_2(u) f^\pm(u-i\h)&=&f^\pm(u) k^\pm_2(u)\ ,\nn\\
k^\pm_1(u) e^\pm(u)&=&e^\pm(u+i\h) k^\pm_1(u)\ ,\nn\\
k^\pm_2(u) e^\pm(u)&=&e^\pm(u-i\h) k^\pm_2(u)\ .\label{com-shift}
\eea
The natural map
$\Delta' L= L{\dot\ot} L $
implies the following universal comultiplication rules
for the Gauss coordinates (see \r{comul-L-univ}):
\bea
\Delta e^\pm(u) &=&
 e^\pm(u')\ot 1 +
 \sum_{p=0} ^{\infty}(-1)^p
  \left(f^\pm(u'-i\tih)\right)^{p} h^\pm(u')\otimes
\left(e^\pm(u'') \right)^{p+1},   \nn\\
\Delta f^\pm(u)&=&
1\otimes f^\pm(u'') +
 \sum_{p=0} ^{\infty} (-1)^p
       \left(f^\pm(u')\right)^{p+1} \otimes
 \tilde h^\pm(u'')\left(e^\pm(u''-i\tih)
\right)^p,   \nn\\
\Delta h^\pm(u)&=&\sum_{p,p'=0}^\infty (-1)^{p+p'}
\left(\left(f^\pm(u')\right)^{p} \ot \left(e^\pm(u''+i\tih) \right)^p
\right)\times \nn\\
&\times&\left(h^\pm(u') \ot h^\pm(u'') \right)
\left(\left(f^\pm(u')\right)^{p'} \ot \left(e^\pm(u''-i\tih) \right)^{p'}
\right)\ .\label{com-gener}
\eea
We see from \r{com-gener} that the choice of the Riemann problem for the
algebra of formal currents enables one to reconstruct also the
comultiplication structure of the algebra.
We will consider this ideology on some examples.
\bigskip

\noindent
{\bf \razdel
The twisting of the Yangian algebras.}
Besides the comultiplication for the Gauss coordinates given by the formulas
\r{com-e-fu}--\r{comul-h} and related to the natural comultiplication
in terms of $L$-operators there is exist another comultiplication
introduced firstly in the paper \cite{Drnew} for the quantized affine
algebras.
\bea
\Delta e(u)&=&e(u')\ot 1   + h^+(u')  \ot e(u'')\ ,   \nn\\
\Delta f(u)&=&1\ot f(u'') + f(u')    \ot h^-(u'')\ ,\nn\\
\Delta h^\pm(u)&=&h^\pm(u') \ot h^\pm(u'')\ .\label{com-Dr}
\eea
This comultiplication describes the coproduct
of the total currents $e(u)$ and $f(u)$ and at the first sight
do not related to $L$-operator formulation of the corresponding
deformed algebra. It was shown in the papers \cite{KT,KT1}
that the corresponding to \r{com-Dr} Hopf algebra
can be obtained as twisted Hopf algebra, which is equivalent to
infinite limit of the shifting automorphism in the quantum
Weyl group. We would like to explain this twisting procedure
on the example of the algebra $\Aya$. The twisting of the
algebra $\DY$ can be considered analogously. The essential
part of the construction will be changing of the Riemann problems
and comultiplication formulas.

Fix the parameter $a\in\RR$.
Consider the following automorphism
\beq
\omega_a\left( e(u)\right)=\ee^{iau} e(u), \quad
\omega_a\left( f(u)\right)=\ee^{-iau} f(u), \quad
\omega_a\left( h^\pm(u)\right)=\ee^{\pm c\h a/2} h^\pm(u),
\label{auto-ful}
\eeq
which obviously conserve the commutation relations
\r{CC1}. The automorphism \r{auto-ful} being translated
to the formal generators $\he_\la$,  $\hf_\la$ and  $\hhh_\la$
 takes the form:
\bea
\omega_a\left( \he_\la\right)&=& \he_{\la+a}, \quad
\omega_a\left( \hf_\la\right)= \hf_{\la-a}, \nn\\
\omega_a\left( \hhh_\la\right)&=&\hhh_\la\left[
\ee^{c\h a/2}\theta(\la)+\ee^{-c\h a/2}\theta(-\la)\right]
+4\h^{-1}\sh(c\h a/2)\delta(\la)\ .
\label{auto-har}
\eea
Appearing of the
$\delta$-function in \r{auto-har} is possible since we consider
the elements of the algebra $\Aya$ as the integrals of the formal
generators with exponentially decreasing weight functions.

The $R$-matrix will change and its value at the critical point
yield the following modification of the commutation relations
\r{com-shift}:
\bea
\ee^{a\h} k^\pm_1(u) f^\pm(u+i\h)&=&f^\pm(u) k^\pm_1(u)\ ,\nn\\
\ee^{-a\h}k^\pm_2(u) f^\pm(u-i\h)&=&f^\pm(u) k^\pm_2(u)\ ,\nn\\
k^\pm_1(u) e^\pm(u)&=&\ee^{-a\h}e^\pm(u+i\h) k^\pm_1(u)\ ,\nn\\
k^\pm_2(u) e^\pm(u)&=&\ee^{a\h} e^\pm(u-i\h) k^\pm_2(u)\ .\label{*4}
\eea
The comultiplication in the twisted algebra due to \r{*4} reads:
\bea
\Delta e^\pm(u) &=&
 e^\pm(u')\ot 1 +
 \sum_{p=0} ^{\infty}(-1)^p
  \left(\ee^{-a\h}f^\pm(u'-i\tih)\right)^{p} h^\pm(u')\otimes
\left(e^\pm(u'') \right)^{p+1},   \nn\\
\Delta f^\pm(u)&=&
1\otimes f^\pm(u'') +
 \sum_{p=0} ^{\infty} (-1)^p
       \left(f^\pm(u')\right)^{p+1} \otimes
h^\pm(u'')\left(\ee^{a\h}e^\pm(u''-i\tih)
\right)^p,   \nn\\
\Delta h^\pm(u)&=&\sum_{p,p'=0}^\infty (-1)^{p+p'}
\left(\left(f^\pm(u')\right)^{p} \ot \left(\ee^{-a\h}e^\pm(u''+i\tih) \right)^p
\right)\times \nn\\
&\times&\left(h^\pm(u') \ot h^\pm(u'') \right)
\left(\left(f^\pm(u')\right)^{p'} \ot \left(\ee^{a\h}
e^\pm(u''-i\tih) \right)^{p'}
\right)\ .\label{com-twist}
\eea

It is obvious that the automorphism \r{auto-ful} changes the
asymptotical  properties of the currents $e^\pm(u)$ and
$f^\pm(u)$. These generating function defined by the Riemann
problems \r{int-f} are analytical in the corresponding domains
of the complex plane $u$ and are decreasing as $u^{-1}$ functions
when $u\to\mp\infty$, respectively.
 The transformed currents
$\omega_a\sk{e^\pm(u)}$ and  $\omega_a\sk{f^\pm(u)}$
will have changed asymptotics
\beq
\omega_a\sk{e^\pm(u)}\sim\ee^{\pm a|{\rm Im}\, u|}, \quad
\omega_a\sk{f^\pm(u)}\sim\ee^{\mp a|{\rm Im}\, u|}, \quad
\Im u\to \mp\infty\ .
\label{asympt}
\eeq
Consider now the limit of the twisted algebra when twisting parameter
$a\to+\infty$. From \r{asympt} it is clear that
$$
\lim{a\to\infty}\omega_a\sk{e^-(u)}=0\,\quad
\lim{a\to\infty}\omega_a\sk{f^+(u)}=0\ .
$$
Because of the Ding-Frenkel relations \r{total} we have
$$
\lim{a\to\infty}e^+(u)=e(u-ic\h/4)\ ,\quad
\lim{a\to\infty}f^-(u)=-f(u-ic\h/4)\ .
$$
So we conclude that the limit of the twisted algebra is an algebra
with the commutation relation given by \r{CC1} and the following
Riemann problem:
\beq
e^-(u)=f^+(u)=0,\quad e(u)=e^+(u+ic\h/4),\quad
f(u)=-f^-(u+ic\h/4).
\label{new-real}
\eeq
From the general formulas \r{com-gener}
the comultiplication \r{com-Dr} follows.
As result we claim that Drinfeld's new realization is the
deformed current algebra with Riemann problem given in
\r{new-real}.

\bigskip

\noindent
{\bf \razdel
The algebra $\Ael$ \cite{KLP3}.}
Let us consider the following generalization of the total current
algebra given by the commutation relations \r{CC1}.
\bea
\H^+(u) \H^-(v)&=&
{
\sh\,\pi\eta(u-v-i\tih(1+c/2))
\over
\sh\,\pi\eta(u-v+i\tih(1-c/2))
}\times \nn\\
&\times&
{
\sh\,\pi\eta'(u-v+i\tih(1+c/2))
\over
\sh\,\pi\eta'(u-v-i\tih(1-c/2))
}
\H^-(v) \H^+(u)\ ,\label{h+h-c}\\
\H^\pm(u) \H^\pm(v)&=&
{
\sh\,\pi\eta(u-v-i\tih)
\sh\,\pi\eta'(u-v+i\tih)
\over
\sh\,\pi\eta(u-v+i\tih)
\sh\,\pi\eta'(u-v-i\tih)
}\
\H^\pm(v) \H^\pm(u)\ ,\label{hhc}\\
\H^\pm(u)\E(v)&=&
{
\sh\,\pi\eta(u-v-i\tih(1\pm c/4))
\over
\sh\,\pi\eta(u-v+i\tih(1\mp c/4))
}\
\E(v)\H^\pm(u)\ , \label{hec}\\
\H^\pm(u)\F(v)&=&
{
\sh\,\pi\eta'(u-v+i\tih(1\pm c/4))
\over
\sh\,\pi\eta'(u-v-i\tih(1\mp c/4))
}\
\F(v) \H^\pm(u)\ , \label{hfc}\\
\E(u)\E(v)&=&
{
\sh\,\pi\eta(u-v-i\tih)
\over
\sh\,\pi\eta(u-v+i\tih)
}\  \E(v)\E(u)\ ,\label{eec}\\
\F(u)\F(v)&=&
{
\sh\,\pi\eta'(u-v+i\tih)
\over
\sh\,\pi\eta'(u-v-i\tih)
}\ \F(v)\F(u)\ ,\label{ffc}\\
{[}\E(u),\F(v){]}&=& {\tih}\left[
\delta\left(u-v+\fract{ic\tih}{2}\right)
\H^+\left(u+\fract{ic\tih}{4}\right)\right.-\nn\\
&-&\left.
\delta\left(u-v-\fract{ic\tih}{2}\right)
\H^-\left(v+\fract{ic\tih}{4}\right)
\right]  ,            \label{efc}
\eea
where the $\delta$-function is defined by the formula
\r{delta-con} and the periods of the trigonometric functions
$\eta$ and $\eta'$ are related as follows:
\beq
{1\over\eta'}-{1\over\eta}=-\h c\ .
\label{eta-rel}
\eeq
Last relation is necessary in order to make the commutation relations
given by \r{h+h-c}--\r{efc}   self consistent.
The ``rational'' commutation relations \r{CC1} can be obtained from
\r{h+h-c}--\r{efc} by the degeneration $\eta\to0$. Note that due
to \r{eta-rel} also $\eta'\to0$ in this limit.
Note that
the formulas \r{h+h-c}--\r{efc} differ from the formulas
given in \cite{KLP3} by reversing the sign of the central element.
The algebra given by the commutation relations \r{h+h-c}--\r{efc}
is associated with simple Lie algebra $\frak{sl}_2$ but can be
formulated for arbitrary simply-laced Lie algebra,
see for example \cite{HZD}.

One can see that the most unusual feature of these relations is the
 presence of two periods $\eta$ and $\eta'$ in trigonometric functions
 playing the role of structure constants. Let us make a comment on
the appearance of such phenomena.
If we put $c=0$  then the periods coincide and the
 relations are analogous to the ones for $c=0$ quantum affine algebra
with $q=\ee^{i\pi\eta\h}$
 in variables $z=\ee^{\pi \eta u}$, $w=\ee^{\pi\eta v}$.
The question is how to input the
 central charge to have nontrivial representation theory. Let us view
 first the classical picture.  In the limit $\h\to 0$  $c=0$ the commutation
 relations \r{h+h-c}--\r{efc} are the relations for Lie algebra  of
 $\frak{sl}_2$ valued (generalized) functions over $z$ vanishing at
 $|\mbox{Re}\ z| \to\pm\infty$ (the Cartan currents tend to $\pm 1$)
\cite{KLPST}
 $$[h^\pm(u),e(v)]=2\pi i\eta\cth \pi\eta(u-v)e(v)\ ,\qquad
   [h^\pm(u),f(v)]=-2\pi i\eta\cth \pi\eta(u-v)e(v)\ ,$$
 \beq
   [e(u),f(v)]=\delta(u-v)(h^+(u)-h^-(v)\ .
\label{717}
\eeq
The application of the Fourier transform to the
generating functions $e(u)=e\ot 2\pi\delta(u-z)$,
$f(u)=f\ot 2\pi\delta(u-z)$, $h^+(u)= h\ot i\pi\eta\cth\pi\eta(u-z)$,
$h^-(u)= h\ot i\pi\eta\cth\pi\eta(u-z-i/\eta)$
which satisfy \r{717}:
$$
e(u)=\int_{-\infty}^{+\infty}d\la\ \ee^{-i\la u}\he_\la\ ,\quad
h^\pm(u)=\int_{-\infty}^{+\infty}d\la\
\ee^{-i\la u}\frac{\hhh_\la}{1-\ee^{\pm\la/\eta}}\ ,\quad
f(u)=\int_{-\infty}^{+\infty}d\la\ \ee^{-i\la u}\hf_\la\ ,$$
turns these commutation relations into the standard form
$$[\hhh_\la,\he_\mu]=2\he_{\la+\mu}\ ,\quad
 [\hhh_\la,\hf_\mu]=-2\hf_{\la+\mu}\ ,\quad
 [\he_\la,\hf_\mu]=\hhh_{\la+\mu} $$
 which admits the standard central extension
$$ [\he_\la,\hf_\mu]=\hhh_{\la+\mu}+\delta(\la+\mu)c\ ,\quad
 [\hhh_\la,\hhh_\mu]=2\delta(\la+\mu)c\ .$$
Let us look to the value of corresponding cocycle on the fields
 $h^+(u_i)$ \cite{KLPST}
\beq
B(h^+(u_1),h^+(u_2))
= 2i\pi\eta^2\left(
{\pi\eta u\over\sh^2\,\pi\eta u}
-\cth\,\pi\eta u\right).
\label{cocycle}
\eeq
The first term in r.h.s. of \r{cocycle} is no longer periodical function
 with a period $1/\eta$. Moreover, the cocycle can be written in term of
the integral over the border of a strip $\Pi: 0<\Im z<1/\eta$ and
 derivatives over the period $1/\eta$:
\beq
B(x\otimes \varphi(z),y\otimes \psi(z))={\eta^2\over 4\pi}
\int_{\partial \Pi}dz \left({d\psi(z)\over d\eta}\varphi(z)-
\psi(z){d\varphi(z)\over d\eta}
\right) \ \langle x,y\rangle\ ,
\label{3.2}
\eeq
where $\langle\ ,\ \rangle$ is Killing form, $x,y\in\frak{sl}_2$.
These arguments signify that the central extension for a quantum
 algebra should be achieved via the finite shift of the periods
 of the trigonometric functions in the defining relations \r{h+h-c}--
\r{efc}.

We attach to the formal current algebra \r{h+h-c}--\r{efc}
the Riemann problem of the following type.

Given meromorphic function $g(u)$, we would like to find two functions
$g^\pm(u)$ satisfying the following conditions:
\begin{itemize}
\item[$(i)$]$g(u)=g^+(u)-g^-(u)$;
\item[$(ii)$]$g^\pm(u)$ are piecewise analytical functions. More
precisely, $g^\pm(u)$ are analytical on the complement to some collection
of the horizontal lines in a complex plane $u$;
\item[$(iii)$]$g^+(u)$ ($g^+(u)$) has a boundary value on the
lower (upper) boundaries of the corresponding strips;
\item[$(iv)$]the following relations hold in any strip of analyticity
of $g^\pm(u)$
$$
g^-(u)=-g^+(u-i\eta)\ ,\quad g^+(u)=-g^-(u+i\eta)\ ,
$$
where $g^+(u-i\eta)$ and $g^-(u+i\eta)$ are the analytical continuations
of $g^+(u)$ and $g^-(u)$.
\end{itemize}

The solution to this Riemann problem is given by the following
integrals over the horizontal lines close to the point $u$:
$$
g^\pm(u)=\pi\eta \int_{\Im v\,\lessgtr\,\Im u}{dv\over 2\pi i}\
{g(v)\over \sh\,\pi\eta(u-v)}\ .
$$
More precisely, the Riemann problems for the
currents $\E(u)$ and $\F(u)$ are chosen with the different periods
$1/\eta$, $1/\eta'$ and certain shifts of the spectral parameter:
\bea
e^\pm(u)&=&\h^{-1}\sin\,\pi\eta\tih
\int_{\Im(u-v)\,\lessgtr\,\pm c\h/4}
{dv\over2\pi i}\
 { \E(v)\over\sh\,\pi\eta
(u-v\mp ic\tih/4)}\ ,\label{2e}\\
f^\pm(u)&=&\h^{-1}\sin\,\pi\eta'\tih \int_{\Im(u-v)\,\lessgtr\,\pm c\h/4}
{dv\over2\pi i}\  { \F(v)\over\sh\,\pi\eta'
(u-v\pm ic\tih/4)}\ .\label{2f}
\eea
The coefficients in front of the integrals in \r{2e} and \r{2f}
are chosen from the technical reasoning and the condition $(i)$
have in this case the form:
\bea
e^+(u+ic\h/4)-e^-(u-ic\h/4)&=&
{\sin\,\pi\eta\h\over\pi\eta\h}\E(u)\ ,\nn\\
f^+(u-ic\h/4)-f^-(u+ic\h/4)&=&
{\sin\,\pi\eta'\h\over\pi\eta'\h}\F(u)\ .\label{EF}
\eea

These Riemann problems are in accordance with the commutation relations
\r{h+h-c}--\r{efc} in the sense that they
yields currents $e^\pm(u)$ and $f^\pm(u)$ which
satisfy the commutation
relations  without the integral terms.
As well as in the case
of the Yangian algebras the Riemann problem is not formulated for
the Cartan  currents and we set
$$
\tilde h^\pm(u)=
{\sin\,\pi\eta'\h\over \pi\eta'\h}\H^\pm(u),\quad
h^\pm(u)=
{\sin\,\pi\eta\h\over \pi\eta\h}H^\pm(u).
$$

The commutation relations between currents $e^+(u)$, $f^+(u)$ and
$h^+(u)$ are given by the relations \r{ef}-\r{hh}.
\bea
e^+
(u_1)f^+(u_2)&-&f^+(u_2)e(u_1)
={\sh\,i\pi\eta'\tih\over\sh\,\pi\eta'u}h^+(u_1)-
{\sh\,i\pi\eta\tih\over\sh\,\pi\eta u}\tilde h^+(u_2),\label{ef}\\
\sh\,\pi\eta(u+i\tih)h^+(u_1)e^+(u_2)&-&
\sh\,\pi\eta(u-i\tih)e^+(u_2)h^+(u_1)=\nn\\
&=&\sh(i\pi\eta\tih)
\{h^+(u_1),e^+(u_1)\},\label{he}\\
\sh\,\pi\eta'(u-i\tih)h^+(u_1)f^+(u_2)&-&
\sh\,\pi\eta'(u+i\tih)f^+(u_2)h^+(u_1)=\nn\\
&=&-\sh(i\pi\eta'\tih)
\{h^+(u_1),f^+(u_1)\},\label{hf}\\
\sh\,\pi\eta(u+i\tih)e^+(u_1)e^+(u_2)&-&
\sh\,\pi\eta(u-i\tih)e^+(u_2)e^+(u_1)=\nn\\
&=&\sh(i\pi\eta\tih)
\left(e^+(u_1)^2+e^+(u_2)^2\right),\label{ee-univ}\\
\sh\,\pi\eta'(u-i\tih)f^+(u_1)f^+(u_2)&-&
\sh\,\pi\eta'(u+i\tih)f^+(u_2)f^+(u_1)=\nn\\
&=&-\sh(i\pi\eta'\tih)
\left(f^+(u_1)^2+f^+(u_2)^2\right),\label{ff-univ}\\
h^+(u_1)h^+(u_2)&=&
{\sh\,\pi\eta'(u+i\tih)\sh\,\pi\eta(u-i\tih)
\over \sh\,\pi\eta(u+i\tih)\sh\,\pi\eta'(u-i\tih)}
h^+(u_2)h^+(u_1)\ .\label{hh}
\eea
The rest of the relations follows from the condition $(iv)$
for the analytical continuations  which in this case reads as:
$$
e^-(u)=-e^+(u-i/\eta'')\ ,\quad
f^-(u)=-f^+(u-i/\eta'')\ ,\quad
\eta''={2\eta\eta'\over\eta+\eta'}\ .
$$
We also impose analogous relation for the
Cartan currents:
$$
h^-(u)=h^+(u-i/\eta'')\ .
$$

The integral formulas \r{2e} and \r{2f} dictate the following presentation
of the currents $e^\pm(u)$ and $f^\pm(u)$:
\bea
e^\pm(u)&=&\pm{\sin\,\pi\eta\tih\over\pi\eta}
\int_{-\infty}^\infty d\la\
\ee^{-i\la u}\
{\he_\la \ee^{\mp c\tih\la/4}\over 1+\ee^{\mp\la/\eta}}\ ,
\label{Lapl-ee}
\\
f^\pm(u)&=&\pm{\sin\,\pi\eta'\tih\over\pi\eta'}
\int_{-\infty}^\infty d\la\  \ee^{-i\la u}\
{\hf_\la \ee^{\pm c\tih\la/4}
\over 1+\ee^{\mp\la/\eta'}}\ ,
\label{Lapl-ff}
\eea
and due to the formulas \r{EF}
the total currents are given by the Fourier transform:
$$
\E(u)=\h \int_{-\infty}^\infty d\la\ \ee^{-i\la u}\he_\la\ ,\quad
\F(u)=\h \int_{-\infty}^\infty d\la\ \ee^{-i\la u}\hf_\la\ .
$$
Note that in the limit $\eta\to0$ formulas \r{Lapl-ee} and \r{Lapl-ff}
go to \r{Lapl-e+f}.
It is natural to guess
that the Cartan generating functions are given also by the
Fourier integrals of the generators $\hk_\la$:
$$%\beq
h^\pm(u)
={\sin\,\pi\eta\tih\over2\pi\eta}
\int_{-\infty}^\infty d\la\  \ee^{-i\la u}\
\hk_\la \ee^{\pm\la/2\eta''}\ .
%\label{Lapl-h}
$$%\eeq

The generators $\he_\la$ and $\hf_\la$ can be defined by means of the
inversion of the Fourier integrals in \r{Lapl-ee}, \r{Lapl-ff} for
the currents $e^+(u)$ and $f^+(u)$ respectively. This can be done
in the domains of the analyticity of the corresponding currents.
These domains depend on the choice of the representation of the
algebra $\Ael$. For instance, the spectral parameter $z$ of the
evaluation representation shifts these domains.
To obtain the commutation relations between the formal generators
$\he_\la$, $\hf_\la$ and $\hk_\la$ one can use the strip
$\left\{-{1/\eta}-{\tih c/4}<\Im u<-{\tih c/4}\right\}$
as it was done in \cite{KLP3}.
As well as in the case of the algebra $\Aya$ the elements of the
algebra $\Ael$ are the formal integrals similar to \r{425}
(see details in \cite{KLP3}).

Considering the generating functions
$e^\pm(u)$,
$f^\pm(u)$,
$h^\pm(u)=k_1^\pm(u)\left(k_2^\pm(u)\right)^{-1}$
as the Gauss coordinates of the $L$-operators  and using
\r{com-shift} we can obtain from the natural
comultiplication in terms of $L$-operator $\Delta' L(u)=
L(u)\dot\ot L(u)$ the comultiplication map
for these generating functions. The verification of the accordance
of these comultiplication rules with the commutation relations
\r{ef}--\r{hh} bring us to the following phenomena. Let us demonstrate
it on the first term in the
comultiplication formulas for the Cartan generating
functions $h^\pm(u)$ \r{com-gener}:
$\Delta h^\pm(u)=h^\pm(u)\ot h^\pm(u)$.
Other terms can be considered analogously.

The commutation relation \r{hh} can be rewritten in the form
\beq
g(u_1-u_2,\xi-\h c)h^\pm(u_1,\xi)  h^\pm(u_2,\xi)  =
h^\pm(u_1,\xi)  h^\pm(u_2,\xi)
g(u_1-u_2,\xi)\ ,
\label{eq1}
\eeq
where
$$
g(u,\xi)={\sh\,\pi\eta(u-i\h)\over \sh\,\pi\eta(u+i\h)},
\quad \xi=\eta^{-1}.
$$
Due to the fact that $\Delta c= c\ot1+1\ot c=c^{(1)}
+c^{(2)}$
we conclude that the commutation relation
$$
g(u_1-u_2,\xi-\h (c^{(1)}+c^{(2)}))
\Delta h^\pm(u_1,\xi)  \Delta h^\pm(u_2,\xi)  =
\Delta h^\pm(u_1,\xi)  \Delta h^\pm(u_2,\xi)
g(u_1-u_2,\xi)\
$$
will follows from \r{eq1} if and only if the comultiplication
for the generating function is defined as follows:
\beq
\Delta h^\pm(u,\xi)=h^\pm(u,\xi-\h c^{(2)})\ot h^\pm(u,\xi)+\cdots
\label{cccc}
\eeq
By the dots we denoted other terms in the series \r{com-gener}.
Actually, analogous phenomena was observed in the paper
\cite{L} in the theory of quantum Sine-Gordon model and implicitly
used for the definition of the intertwining operators in the
elliptic algebra
$\Apq$ which was served as the dynamical symmetry algebra
of the eight-vertex model
\cite{FIJKMY}.
 The algebra $\Ael$ which is a scaling limit of $\Apq$
was
investigated in \cite{KLP3}.
In the $L$-operator formalism the algebra $\Apq$ as well as $\Ael$
is given by the commutation relation:
$$
R^+(u_1-u_2,\xi-\h c)L(u_1,\eta)L_2(u_2,\eta)=
L_2(u_2,\eta)L_1(u_1,\eta) R^+(u_1-u_2,\xi)\ .
$$
Although the comultiplication map moves the algebra into
the tensor product of the algebra with shifted by the central
element parameter this unusual Hopf structure allows to
construct the intertwining operators for the infinite-dimensional
representations of the algebra $\Ael$ and interpret them as
Zamolodchikov-Faddeev operators in the quantum Sine-Gordon model
\cite{KLP3}.

To conclude we would like to mention the papers \cite{JKOS} and
\cite{ABRR} where the quasi-Hopf twisting of the quantum groups have
been considered. These papers are the realization of the idea
which belongs to C.~Fr\o nsdal \cite{F}
that the elliptic deformations of the
quantum groups can be obtained by the twisting of the corresponding
Hopf algebras and belong to the category of Drinfeld's quasi-Hopf
algebras \cite{DRlen}. In particular, it means that the algebra
$\Ael$ can be obtained by quasi-Hopf
twisting of the algebra $\Aya$ and also explains the comultiplication
formulas like \r{cccc} first considered in \cite{KLP3}. Also the
result of these papers can be considered as a proof that all elliptic
deformations of the Hopf algebras are actually isomorphic for the
generic values of the elliptic parameter, the fact which was
observed in \cite{KLPST} on the level of the classical limit
of the algebra $\Ael$.

\section{Acknowledgement}

S.~Khoroshkin and S.~Pakuliak are grateful for the warm hospitality
to the University Roma~1 where the part of this work was done.
This visit was arranged in the framework of INFN--ITEP and INFN--JINR
exchange programs. S.~Khoroshkin also acknowledge the hospitality
of Max Planck Institut f\"ur Mathematik.
D.~Lebedev is grateful to the Institute GIRARD DESARGUES, University
Claude Bernard Lyon 1 for warm hospitality where the part of this
work was done.

The work was supported in part by grants
RFBR-96-01-00814a (S.~Kho\-rosh\-kin),
RFBR-95-01-01101 (D.~Lebedev),
RFBR-95-01-01106 (S.~Pakuliak),
INTAS-93-0166-Ext (S.~Khoroshkin, D.~Lebedev),
French MENESRIP/DAEIF (D.~Lebedev),
No. 96-15-96455  for support of scientific schools
(S.~Khoroshkin, D.~Lebedev and S.~Pakuliak)
and by Award  No. RM2-150 of the U.S. Civilian Research \& Development
Foundation (CRDF) for the Independent States of the Former Soviet Union
(S.~Khoroshkin, D.~Lebedev, S.~Pakuliak).

\end{document}